# Asymptotic results for sample autocovariance functions and extremes of integrated generalized Ornstein–Uhlenbeck processes

VICKY FASEN

*Center for Mathematical Sciences, Technische Universität München, D-85747 Garching, Germany. E-mail: fasen@ma.tum.de; url: www-m4.ma.tum.de*

We consider a positive stationary generalized Ornstein–Uhlenbeck process

$$V_t = e^{-\xi_t}\left(\int_0^t e^{\xi_{s-}}\,d\eta_s + V_0\right) \qquad \text{for } t \geq 0,$$

and the increments of the integrated generalized Ornstein–Uhlenbeck process $I_k = \int_{k-1}^k \sqrt{V_{t-}}\,dL_t$, $k \in \mathbb{N}$, where $(\xi_t, \eta_t, L_t)_{t\geq 0}$ is a three-dimensional Lévy process independent of the starting random variable $V_0$. The genOU model is a continuous-time version of a stochastic recurrence equation. Hence, our models include, in particular, continuous-time versions of ARCH(1) and GARCH(1, 1) processes. In this paper we investigate the asymptotic behavior of extremes and the sample autocovariance function of $(V_t)_{t\geq 0}$ and $(I_k)_{k\in\mathbb{N}}$. Furthermore, we present a central limit result for $(I_k)_{k\in\mathbb{N}}$. Regular variation and point process convergence play a crucial role in establishing the statistics of $(V_t)_{t\geq 0}$ and $(I_k)_{k\in\mathbb{N}}$. The theory can be applied to the COGARCH(1, 1) and the Nelson diffusion model.

*Keywords:* continuous-time GARCH process; extreme value theory; generalized Ornstein–Uhlenbeck process; integrated generalized Ornstein–Uhlenbeck process; mixing; point process; regular variation; sample autocovariance function; stochastic recurrence equation

## 1. Introduction

In this paper we develop limit results for stationary positive *generalized Ornstein–Uhlenbeck* (genOU) processes

$$V_t = e^{-\xi_t}\left(\int_0^t e^{\xi_{s-}}\,d\eta_s + V_0\right) \qquad \text{for } t > 0, \tag{1.1}$$

and *integrated genOU* (IgenOU) processes

$$I_t^* = \int_0^t \sqrt{V_{s-}}\,\mathrm{d}L_s \qquad \text{for } t \geq 0, \tag{1.2}$$

where $(\xi_t, \eta_t, L_t)_{t\geq 0}$ is a three-dimensional Lévy process independent of the starting random variable $V_0$, $(\eta_t)_{t\geq 0}$ is a subordinator and $(-L_t)_{t\geq 0}$ is not a subordinator. Here and in general $\int_a^b$ means the integral over $(a, b]$. A three-dimensional Lévy process is characterized by its *Lévy–Khinchine representation* $\mathbb{E}(\mathrm{e}^{\mathrm{i}\langle\Theta,(\xi_t,\eta_t,L_t)\rangle}) = \exp(-t\Psi(\Theta))$ for $\Theta \in \mathbb{R}^3$, where

$$\begin{aligned}\Psi(\Theta) = &-\mathrm{i}\langle\gamma,\Theta\rangle + \frac{1}{2}\langle\Theta, \Sigma\Theta\rangle \\ &+ \int_{\mathbb{R}^3} (1 - \mathrm{e}^{\mathrm{i}\langle\Theta,(x,y,z)\rangle} + \mathrm{i}\mathbf{1}_{\{\langle(x,y,z),(x,y,z)\rangle\leq 1\}}\langle(x,y,z),\Theta\rangle)\,\mathrm{d}\Pi_{\xi,\eta,L}(x,y,z)\end{aligned}$$

with $\gamma \in \mathbb{R}^3$, $\Sigma$ a non-negative definite matrix in $\mathbb{R}^{3\times 3}$ and $\Pi_{\xi,\eta,L}$ a measure on $\mathbb{R}^3$, called *Lévy measure*, which satisfies $\int_{\mathbb{R}^3} \min\{x^2 + y^2 + z^2, 1\}\,\mathrm{d}\Pi_{\xi,\eta,L}(x,y,z) < \infty$ and $\Pi_{\xi,\eta,L}((0,0,0)) = 0$. Further, $\langle\cdot,\cdot\rangle$ denotes the inner product in $\mathbb{R}^3$. A subordinator is a positive Lévy process; we refer to the monographs of Sato (1999) and Applebaum (2004) for more details on Lévy processes. A fundamental contribution to the probabilistic properties of genOU processes is the recent paper of Lindner and Maller (2005).

GenOU processes are applied in various areas, for example, in financial and insurance mathematics, or mathematical physics; we refer to Carmona *et al.* (1997, 2001), and Donati-Martin *et al.* (2001) for an overview of applications. Processes of this class are used as stochastic volatility models in finance (cf. Barndorff-Nielsen and Shephard (2001)) and as risk models in insurance (cf. Hipp and Plum (2003); Paulsen (2002); Kostadinova (2007)). Continuous-time processes are particularly appropriate models for irregularly-spaced and high-frequency data. A genOU process is a continuous-time version of a stochastic recurrence equation; see de Haan and Karandikar (1989). Practical applications of stochastic recurrence equations are given in Diaconis and Freedman (1999). This means the ARCH(1) process, as solution of a stochastic recurrence equation, can be interpreted as a discrete-time version of a genOU process. A typical example of an IgenOU process is the continuous-time GARCH(1,1) (COGARCH(1,1)) process introduced by Klüppelberg *et al.* (2004) (cf. Example 2.4). On the other hand, Nelson (1990) suggested the approximation of a diffusion by GARCH(1,1) models (cf. Example 2.3). The diffusion model is again an IgenOU process and its volatility process is a genOU process.

We investigate the asymptotic behavior of extremes and the sample autocovariance function, respectively, of

$$H_k = \sup_{(k-1)h \leq t \leq kh} V_t \qquad \text{for } k \in \mathbb{N} \tag{1.3}$$

and some $h > 0$, of $(V_t)_{t \geq 0}$ and of the stationary increments

$$I_k = I_{kh}^* - I_{(k-1)h}^* = \int_{(k-1)h}^{kh} \sqrt{V_{t-}}\, \mathrm{d}L_t \qquad \text{for } k \in \mathbb{N}, \tag{1.4}$$

of $(I_t^*)_{t \geq 0}$. Including continuous-time versions of ARCH(1) and GARCH(1, 1) processes, we derive similar results, as in Davis and Mikosch (1998) and Mikosch and Stărică (2000), who investigated the asymptotic behavior of extremes and the sample autocovariance functions of ARCH(1) and GARCH(1, 1) processes, for $(V_t)_{t \geq 0}$ and $(I_k)_{k \in \mathbb{N}}$.

In this paper we present only theoretical results. One reason is that financial time series often have finite variance but infinite fourth moment. If the IgenOU or the genOU processes have these properties, then the normalized sample autocovariance functions of $(I_k)_{k \in \mathbb{N}}$ and $(V_t)_{t \geq 0}$, respectively, converge to an infinite variance stable distribution (see Section 4.3). The structure of these stable distributions is complex, and it is not clear how to compute them analytically. Hence, it is also difficult to calculate any confidence intervals from these results. Further, we restrict our attention to only qualitative results, since the inference, estimation and testing of a genOU and an IgenOU process is not fully developed. First steps in estimation procedures of the COGARCH(1, 1) process are given in Haug *et al.* (2007), Maller *et al.* (2008) and Müller (2007).

The paper is organized as follows: We start, in Section 2, with a detailed analysis of the genOU and the IgenOU model used in this paper. This analysis includes sufficient conditions for model assumptions and examples. The regular variation of these processes, stated in Section 3.1, is crucial to proving the convergence of relevant point processes. These conclusions agree with the empirical findings of heavy tailed logarithmic returns of financial time series. Section 3.2 concerns mixing properties of $(V_t)_{t \geq 0}$, $(H_k)_{k \in \mathbb{N}}$ and $(I_k)_{k \in \mathbb{N}}$.

First, we derive in Section 4 the convergence of point processes based on $(H_k)_{k \in \mathbb{N}}$ and $(I_k)_{k \in \mathbb{N}}$. These results we use to develop the extremal behavior of $(H_k)_{k \in \mathbb{N}}$ in Section 4.1, the asymptotic behavior of $(I_t^*)_{t \geq 0}$, in the form of a central limit result, in Section 4.2, and the asymptotic behavior of the sample autocovariance functions of $(V_t)_{t \geq 0}$ and $(I_k)_{k \in \mathbb{N}}$ in Section 4.3. One important conclusion is that $(H_k)_{k \in \mathbb{N}}$ and $(I_k)_{k \in \mathbb{N}}$ exhibit extremal clusters, which are often observed in financial time series. Finally, the proofs of the results are included in Appendices A and B.

We shall use the following standard notations: $\mathbb{R}_+ = (0, \infty)$. For real functions $g$ and $h$ we abbreviate $g(t) \sim h(t)$ for $t \to \infty$, if $g(t)/h(t) \to 1$ for $t \to \infty$. For $x \in \mathbb{R}$ we set $x^+ = \max\{x, 0\}$ and $x^- = \max\{0, -x\}$. For a vector $\mathbf{x} \in \mathbb{R}^k$ we also denote by $|\mathbf{x}|_\infty = \max\{|x_1|, \ldots, |x_k|\}$ the maximum norm. We write $X \stackrel{d}{=} Y$, if the distributions of the random variables $X$ and $Y$ coincide. Provided that $\mathbb{E}(\mathrm{e}^{-v\xi_1})$ is finite for $v > 0$ we set

$$\Psi_\xi(v) = \log \mathbb{E}(\mathrm{e}^{-v\xi_1}).$$

Then $\mathbb{E}(\mathrm{e}^{-v\xi_t}) = \mathrm{e}^{t\Psi_\xi(v)}$ is finite for all $t \geq 0$; see Sato (1999), Theorem 25.17.

## 2. Model assumptions and examples

### 2.1. Model assumptions

Throughout the paper we assume that the genOU process satisfies at least condition (A) as below.

***Condition (A).*** *The stochastic process $(V_t)_{t\geq 0}$ is a stationary positive càdlàg version of the genOU process in (1.1). Further, the stationary distribution $V_0$ has a Pareto-like tail with index $\alpha > 0$, that is, $\mathbb{P}(V_0 > x) \sim C x^{-\alpha}$ as $x \to \infty$ for some $C > 0$.*

This is a natural condition; see Proposition 2.1 below for a precise formulation of sufficient assumptions. We will assume either condition (B) or (C) hereafter depending on whether we investigate probabilistic properties of the genOU process or the IgenOU process.

***Condition (B).*** *There exist $\alpha > 0$ and $d > \alpha$ such that*

$$\Psi_\xi(\alpha) = 0 \quad and \quad \Psi_\xi(d) < \infty. \tag{2.1}$$

*Furthermore, for some $h > 0$,*

$$\mathbb{E}\left| e^{-\xi_h} \int_0^h e^{\xi_{s-}} \, d\eta_s \right|^d < \infty. \tag{2.2}$$

Condition (B) stems from the application of results for stochastic recurrence equations of Kesten (1973) and Goldie (1991) to the equation $V_{(k+1)h} = A_{(k+1)h}^{kh} V_{kh} + B_{(k+1)h}^{kh}$ for $k \in \mathbb{N}$, where

$$A_t^s = e^{-(\xi_t - \xi_s)} \quad \text{and} \quad B_t^s = e^{-\xi_t} \int_s^t e^{\xi_{u-}} \, d\eta_u \qquad \text{for } 0 \leq s < t.$$

A conclusion of de Haan and Karandikar (1989) (see also Carmona *et al.* (1997)) is that $(V_t)_{t\geq 0}$ is a time-homogenous Markov process and $(A_{(k+1)h}^{kh}, B_{(k+1)h}^{kh})_{k\in\mathbb{N}}$ is an i.i.d. sequence.

***Condition (C).*** *Suppose condition (B) is satisfied. Furthermore, for $k \in \mathbb{N}$, $\mathbb{E}|L_1| < \infty$,*

$$\mathbb{E}\left| \int_0^h e^{-\xi_{t-}/2} \, dL_t \right|^{2\max\{1,d\}} < \infty,$$

$$\mathbb{E}\left| \int_{(k-1)h}^{kh} e^{-\xi_{t-}/2} \left( \int_0^{t-} e^{\xi_{s-}} \, d\eta_s \right)^{1/2} dL_t \right|^{2\max\{1,d\}} < \infty. \tag{2.3}$$

This condition arises from the following decomposition of

$$I_k = \int_{(k-1)h}^{kh} \sqrt{A_{t-}^{(k-1)h} V_{(k-1)h} + B_{t-}^{(k-1)h}} \, \mathrm{d}L_t.$$

Thus, assumption (2.3) is equivalent to

$$\mathbb{E}\left| \int_{(k-1)h}^{kh} \sqrt{A_{t-}^{(k-1)h}} \, \mathrm{d}L_t \right|^{2\max\{1,d\}} < \infty \quad \text{and} \quad \mathbb{E}\left| \int_{(k-1)h}^{kh} \sqrt{B_{t-}^{(k-1)h}} \, \mathrm{d}L_t \right|^{2\max\{1,d\}} < \infty.$$

Theorem 4.5 of Lindner and Maller (2005) presents sufficient conditions for (A), which is included in the next proposition.

**Proposition 2.1.** *Let $(V_t)_{t\geq 0}$ be the genOU process in (1.1). When $\xi$ is of finite variation, we assume additionally that the drift of $\xi$ is non-zero, or that there is no $r > 0$ such that the support of the Lévy measure of $\xi$ is concentrated on $r\mathbb{Z}$.*

(a) *Suppose there exist $\alpha > 0$, $d > \alpha$, $p,q > 1$ with $1/p + 1/q = 1$ such that*

$$\Psi_\xi(\alpha) = 0, \qquad \mathbb{E}(\mathrm{e}^{-\max\{1,d\}p\xi_1}) < \infty \quad \text{and} \quad \mathbb{E}|\eta_1|^{q\max\{1,d\}} < \infty. \tag{2.4}$$

*Then there exists a version of $V$ satisfying condition (A), and (B) holds.*

(b) *Suppose there exist $p_i, q_i > 1$ with $1/p_i + 1/q_i = 1$, $i = 1,2$, $\alpha > 0$ and $d > \alpha$ such that*

$$\begin{aligned}\Psi_\xi(\alpha) &= 0, \qquad \mathbb{E}(\mathrm{e}^{-p_1 p_2 \max\{1,d\}\xi_1}) < \infty, \\ \mathbb{E}|\eta_1|^{q_1 p_2 \max\{1,d\}} &< \infty, \qquad \mathbb{E}|L_1|^{2q_2 \max\{1,d\}} < \infty.\end{aligned} \tag{2.5}$$

*Then there exists a version of $V$ satisfying condition (A), and (C) holds.*

### 2.2. Examples

*Example 2.2 (Ornstein–Uhlenbeck process).* The Lévy-driven Ornstein–Uhlenbeck process

$$V_t = \mathrm{e}^{-\lambda t}\left(\int_0^t \mathrm{e}^{\lambda s}\,\mathrm{d}\eta_s + V_0\right) \qquad \text{for } t \geq 0,$$

is a simple example of a genOU process. Since $\Psi_\xi(s) = -s\lambda < 0$ for $s > 0$ the assumption (2.1) cannot be satisfied. Hence, this process is not included in the framework of this paper; we refer to Fasen *et al.* (2006) for more details on extreme value theory of Lévy-driven Ornstein–Uhlenbeck processes.

***Example 2.3 (Nelson's diffusion model).*** In the diffusion model of Nelson (1990) the volatility process is the stationary solution of the SDE

$$\mathrm{d}V_t = \lambda(a - V_t)\,\mathrm{d}t + \sigma V_t\,\mathrm{d}W_t^{(1)} \qquad \text{for } t \geq 0, \tag{2.6}$$

where $\lambda, a, \sigma > 0$ and $(W_t^{(1)})_{t\geq 0}$ is a Brownian motion. Then Nelson (1990) models logarithmic asset prices of financial time series by

$$I_t^* = \int_0^t \sqrt{V_t}\,\mathrm{d}W_t^{(2)} \qquad \text{for } t \geq 0,$$

where $(W_t^{(2)})_{t\geq 0}$ is a Brownian motion independent of $(W_t^{(1)})_{t\geq 0}$. Theorem 52 in Protter (2004), page 328, gives that $(V_t)_{t\geq 0}$ is a genOU process with representation

$$V_t = \mathrm{e}^{-\xi_t}\left(\lambda a \int_0^t \mathrm{e}^{\xi_s}\,\mathrm{d}s + V_0\right) \qquad \text{for } t \geq 0,$$

where $\xi_t = -\sigma W_t^{(1)} + (\sigma^2/2 + \lambda)t$. In this case

$$(\xi_t, \eta_t, L_t) = (-\sigma W_t^{(1)} + (\sigma^2/2 + \lambda)t, \lambda at, W_t^{(2)}).$$

Here, we do not take left limits of $\xi$ (or $V$, resp.) in the representation of the genOU process, since $\xi$ has continuous sample paths from the Brownian motion. Furthermore,

$$\Psi_\xi(v) = -\left(\frac{1}{2}\sigma^2 + \lambda\right)v + \frac{\sigma^2}{2}v^2 \qquad \text{for } v \in \mathbb{R},$$

so that for $\alpha = 1 + 2\lambda/\sigma^2$ we have $\Psi_\xi(\alpha) = 0$. Hence, there exists a version of $V$ and $I$ satisfying assumptions (A)–(C) for any $d > \alpha$, $p_i, q_i > 1$ with $1/p_i + 1/q_i = 1$, $i = 1, 2$.

***Example 2.4 (*COGARCH$(1,1)$ *model).*** Let $(\xi_t)_{t\geq 0}$ be a spectrally negative Lévy process with representation

$$\xi_t = ct - \sum_{0<s\leq t} \log(1 + \lambda \mathrm{e}^c (\Delta L_s)^2) \qquad \text{for } t \geq 0,$$

where $c > 0$, $\lambda \geq 0$ and $(L_t)_{t\geq 0}$ is a Lévy process. Then the volatility process of the COGARCH$(1,1)$ process as defined in Klüppelberg *et al.* (2004) (we use only the right-continuous version) is given by

$$V_t = \mathrm{e}^{-\xi_t}\left(\beta \int_0^t \mathrm{e}^{\xi_{s-}}\,\mathrm{d}s + V_0\right) \qquad \text{for } t \geq 0 \tag{2.7}$$

and $\beta > 0$. With this definition the COGARCH$(1,1)$ process has the representation

$$I_t^* = \int_0^t \sqrt{V_{t-}}\,\mathrm{d}L_t \qquad \text{for } t \geq 0.$$

In contrast to the Nelson diffusion model, $\xi$ and $L$ are dependent here.

(a) If there exists an $\alpha > 0$ and $d > \alpha$ such that

$$\Psi_\xi(\alpha) = 0 \quad \text{and} \quad \mathbb{E}|L_1|^{2d} < \infty, \tag{2.8}$$

then a stationary version of $V$ exists, whose marginal distribution is regularly varying with index $\alpha$; see Klüppelberg *et al.* (2006), Theorem 6. Hence, (A) and also (B) follow.

(b) If we assume that there exist an $\alpha > 0$ and some $d > \alpha$ such that

$$\Psi_\xi(\alpha) = 0 \quad \text{and} \quad \mathbb{E}|L_1|^{\max\{4d,1\}} < \infty, \tag{2.9}$$

then additionally (C) is satisfied.

## 3. Preliminary results

### 3.1. Regular variation

The tail behavior of the stationary distribution has a crucial impact on the extremes of a stationary process. But the dependence of large values in successive variables also has an influence on the extremal behavior of stochastic processes. A possible model for large values in different components is, in our case, regular variation of the continuous-time process $V$. Regular variation of stochastic processes was studied by de Haan and Lin (2001) and Hult and Lindskog (2007). Before we present the definition we require some notation. Let $\mathbb{D}$ be the space of all càdlàg functions on $[0,1]$ equipped with the $J_1$-metric which gives the Skorokhod topology (cf. Billingsley (1999)) and $\mathbb{S}_\mathbb{D} = \{\mathbf{x} \in \mathbb{D} : |\mathbf{x}|_\infty = 1\}$ is the unit sphere in $\mathbb{D}$ equipped with the subspace topology, where $|\mathbf{x}|_\infty = \sup_{0 \le t \le 1} |x_t|$. The symbol $\mathcal{B}$ denotes the Borel $\sigma$-algebra and $\stackrel{u \to \infty}{\Longrightarrow}$ weak convergence as $u \to \infty$.

**Definition 3.1.** *A stochastic process $\mathbf{X} = (X_t)_{0 \le t \le 1}$ with sample paths in $\mathbb{D}$ is said to be* regularly varying *with index $\alpha > 0$, if there exists a probability measure $\sigma$ on $\mathcal{B}(\mathbb{S}_\mathbb{D})$ such that for every $x > 0$,*

$$\frac{\mathbb{P}(|\mathbf{X}|_\infty > ux, \mathbf{X}/|\mathbf{X}|_\infty \in \cdot)}{\mathbb{P}(|\mathbf{X}|_\infty > u)} \stackrel{u \to \infty}{\Longrightarrow} x^{-\alpha} \sigma(\cdot) \qquad on\ \mathcal{B}(\mathbb{S}_\mathbb{D}).$$

In this section we consider the tail behavior of $(H_k)_{k \in \mathbb{N}}$ and $(I_k)_{k \in \mathbb{N}}$ described by multivariate regular variation, and we will use these results to derive the convergence of point processes based on these sequences in Section 4. More details and properties on multivariate regularly varying random vectors can be found, for example, in Jessen and Mikosch (2006) and Resnick (1987, 2007).

**Definition 3.2.** *A random vector $\mathbf{X} = (X_1, \ldots, X_k)$ on $\mathbb{R}^k$ is said to be* regularly varying *with index $\alpha > 0$, if there exists a random vector $\Theta$ with values on the unit sphere $\mathbb{S}_{k-1} =$*

$\{\mathbf{x} \in \mathbb{R}^k : |\mathbf{x}|_\infty = 1\}$ such that for every $x > 0$,

$$\frac{\mathbb{P}(|\mathbf{X}|_\infty > ux, \mathbf{X}/|\mathbf{X}|_\infty \in \cdot)}{\mathbb{P}(|\mathbf{X}|_\infty > u)} \overset{u \to \infty}{\Longrightarrow} x^{-\alpha} \mathbb{P}(\Theta \in \cdot) \qquad \text{on } \mathcal{B}(\mathbb{S}_{k-1}).$$

The next theorem shows that the regular variation of $V_0$ has consequences on the processes $(V_t)_{t \geq 0}$, $(H_k)_{k \in \mathbb{N}}$ and $(I_k)_{k \in \mathbb{N}}$.

**Theorem 3.3 (Regular variation).** *Let $(V_t)_{t \geq 0}$ be a genOU process satisfying (A). Further, let $(H_k)_{k \in \mathbb{N}}$ and $(I_k)_{k \in \mathbb{N}}$, respectively, be the stationary processes in (1.3) and (1.4).*

(a) *Suppose (B) is satisfied. Let $\mathbf{V} = (V_t)_{0 \leq t \leq 1}$. Then for every $x > 0$,*

$$\frac{\mathbb{P}(|\mathbf{V}|_\infty > ux, \mathbf{V}/|\mathbf{V}|_\infty \in \cdot)}{\mathbb{P}(|\mathbf{V}|_\infty > u)} \overset{u \to \infty}{\Longrightarrow} x^{-\alpha} \frac{\mathbb{E}(|\mathbf{U}|_\infty^\alpha \mathbf{1}\{\mathbf{U}/|\mathbf{U}|_\infty \in \cdot\})}{\mathbb{E}|\mathbf{U}|_\infty^\alpha} \qquad \text{on } \mathcal{B}(\mathbb{S}_\mathbb{D}),$$

*where $\mathbf{U} = (e^{-\xi_t})_{0 \leq t \leq 1}$.*

(b) *Suppose (B) is satisfied. Let $\mathbf{H}_k = (H_1, \ldots, H_k)$ for $k \in \mathbb{N}$. Then for every $x > 0$,*

$$\frac{\mathbb{P}(|\mathbf{H}_k|_\infty > ux, \mathbf{H}_k/|\mathbf{H}_k|_\infty \in \cdot)}{\mathbb{P}(|\mathbf{H}_k|_\infty > u)} \overset{u \to \infty}{\Longrightarrow} x^{-\alpha} \frac{\mathbb{E}(|\mathbf{m}_k|_\infty^\alpha \mathbf{1}\{\mathbf{m}_k/|\mathbf{m}_k|_\infty \in \cdot\})}{\mathbb{E}|\mathbf{m}_k|_\infty^\alpha} \qquad \text{on } \mathcal{B}(\mathbb{S}_{k-1}),$$

*where*

$$\mathbf{m}_k = \left( \sup_{0 \leq t \leq h} e^{-\xi_t}, \ldots, \sup_{(k-1)h \leq t \leq kh} e^{-\xi_t} \right).$$

*Furthermore,*

$$\mathbb{P}(H_1 > x) \sim \mathbb{E}\left( \sup_{0 \leq s \leq h} e^{-\alpha \xi_s} \right) \mathbb{P}(V_0 > x) \qquad \text{as } x \to \infty.$$

(c) *Suppose (C) is satisfied. Let $\mathbf{I}_k = (I_1, \ldots, I_k)$ for $k \in \mathbb{N}$. Then for every $x > 0$,*

$$\frac{\mathbb{P}(|\mathbf{I}_k|_\infty > ux, \mathbf{I}_k/|\mathbf{I}_k|_\infty \in \cdot)}{\mathbb{P}(|\mathbf{I}_k|_\infty > u)} \overset{u \to \infty}{\Longrightarrow} x^{-2\alpha} \frac{\mathbb{E}(|\mathbf{r}_k|_\infty^{2\alpha} \mathbf{1}\{\mathbf{r}_k/|\mathbf{r}_k|_\infty \in \cdot\})}{\mathbb{E}|\mathbf{r}_k|_\infty^{2\alpha}} \qquad \text{on } \mathcal{B}(\mathbb{S}_{k-1}),$$

*where*

$$\mathbf{r}_k = \left( \int_0^h e^{-\xi_{t-}/2} \, dL_t, \ldots, \int_{(k-1)h}^{kh} e^{-\xi_{t-}/2} \, dL_t \right).$$

*Furthermore,*

$$\mathbb{P}(I_1 > x) \sim \mathbb{E}\left[ \left( \int_0^h e^{-\xi_{t-}/2} \, dL_t \right)^{+} {}^{2\alpha} \right] \mathbb{P}(V_0 > x^2) \qquad \text{as } x \to \infty.$$

## 3.2. Mixing properties

The mixing property of a stochastic process describes the temporal dependence in data. Different kinds of mixing properties have been defined, which are summarized, for example, in the survey paper of Bradley (2005). For the derivation of limit results of point processes in Section 4, one assumption is the asymptotic independence in extrema. Further, mixing is used to prove consistency and asymptotic normality of estimators.

Let $(X_t)_{t\geq 0}$ be a stationary process, $\mathcal{F}_t = \sigma(X_s : s \leq t)$ and $\mathcal{G}_t = \sigma(X_s : s \geq t)$. If

$$\alpha(t) := \sup_{A \in \mathcal{F}_v, B \in \mathcal{G}_{v+t}} |\mathbb{P}(A \cap B) - \mathbb{P}(A)\mathbb{P}(B)| \longrightarrow 0 \qquad \text{as } t \to \infty,$$

then $(X_t)_{t\geq 0}$ is called $\alpha$-mixing. $(X_t)_{t\geq 0}$ is called $\beta$-mixing, if

$$\beta(t) := \sup_{\substack{A_i \in \mathcal{F}_v, i=1,\ldots,I, \\ B_j \in \mathcal{G}_{v+t}, j=1,\ldots,J, \\ \sum_{i=1}^I A_i = \sum_{j=1}^J B_j = \Omega}} \frac{1}{2} \sum_{i=1}^I \sum_{j=1}^J |\mathbb{P}(A_i \cap B_j) - \mathbb{P}(A_i)\mathbb{P}(B_j)| \longrightarrow 0 \qquad \text{as } t \to \infty.$$

The following inequality holds: $2\alpha(t) \leq \beta(t)$. Hence, $\beta$-mixing implies $\alpha$-mixing. $(X_t)_{t\geq 0}$ is called exponentially $\beta$-mixing, if $\beta(t) \leq K e^{-at}$ for some $K, a > 0$ and all $t \geq 0$. Analogous is the definition of exponentially $\alpha$-mixing.

**Proposition 3.4 (Mixing).** *Let $(V_t)_{t\geq 0}$ be a genOU process satisfying* (A) *and* (B). *We assume that $(V_t)$ is simultaneously $\varphi$-irreducible (for some $\sigma$-finite measure $\varphi$). Further, let $(H_k)_{k \in \mathbb{N}}$ and $(I_k)_{k \in \mathbb{N}}$, respectively, be the stationary processes in* (1.3) *and* (1.4).

(a) *Then $(V_t)_{t\geq 0}$ is exponentially $\beta$-mixing and geometrically ergodic.*
(b) *Then $(H_k)_{k \in \mathbb{N}}$ is exponentially $\beta$-mixing and geometrically ergodic.*
(c) *Suppose $(L_t)_{t\geq 0}$ is a Brownian motion independent of $(\xi_t, \eta_t)_{t\geq 0}$. Then $(I_k)_{k \in \mathbb{N}}$ is exponentially $\beta$-mixing and geometrically ergodic.*

### Example 3.5.

(a) Consider the COGARCH$(1,1)$ model of Example 2.4, which satisfies (2.8). Then $(V_t)_{t\geq 0}$ is simultaneously $\lambda$-irreducible, where $\lambda$ denotes the Lebesgue measure (cf. Paulsen (1998), page 142, and Nyrhinen (2001)). Hence, $(V_t)_{t\geq 0}$ and $(I_k)_{k \in \mathbb{N}}$ are exponentially $\beta$-mixing and geometrically ergodic by Proposition 3.4 and Haug *et al.* (2007), Theorem 3.5.

(b) In the Nelson diffusion model (Example 2.3) $(V_t)_{t\geq 0}$ and $(I_k)_{k \in \mathbb{N}}$ are exponentially $\beta$-mixing and geometrically ergodic; see Genon-Catalot *et al.* (2000).

For the derivation of point process results we need the asymptotic independence in extremes as below. It is particularly satisfied for $\alpha$-mixing and $\beta$-mixing sequences (Basrak (2000), Lemma 3.2.9).

**Condition** $\mathcal{A}(c_n)$. *Let $(\mathbf{Y}_k)_{k\in\mathbb{Z}}$ be a strictly stationary sequence of regularly varying random vectors and $0 < c_n \uparrow \infty$ be a sequence of constants satisfying*

$$\lim_{n\to\infty} n\mathbb{P}(|\mathbf{Y}_1|_\infty > c_n) = C \tag{3.1}$$

*for some $C > 0$. There exists a set of positive integers $(r_n)_{n\in\mathbb{N}}$ such that $r_n \to \infty$, $r_n/n \to 0$ as $n \to \infty$ and*

$$\mathbb{E}\exp\left(-\sum_{j=1}^n f(\mathbf{Y}_j/c_n)\right) - \left[\mathbb{E}\exp\left(-\sum_{j=1}^{r_n} f(\mathbf{Y}_j/c_n)\right)\right]^{\lfloor n/r_n \rfloor} \stackrel{n\to\infty}{\longrightarrow} 0 \quad \text{for all } f \in \mathcal{F}_s,$$

*where $\mathcal{F}_s$ is the collection of bounded non-negative step functions on $\overline{\mathbb{R}}^d \setminus \{\mathbf{0}\}$ with bounded support.*

Thus, $(H_k)_{k\in\mathbb{N}}$ and $(I_k)_{k\in\mathbb{N}}$ of Example 2.3 and 2.4, respectively, satisfy condition $\mathcal{A}$ by Theorem 3.3 and Example 3.5.

## 4. Point process convergence and conclusions

In this section we study the weak convergence of point processes of exceedances associated with $(H_k)_{k\in\mathbb{N}}$ and $(I_k)_{k\in\mathbb{N}}$. Point processes are prominent tools to precisely describe the extremal behavior of stochastic processes (see Resnick (1987, 2007)). They can be used to determine the limit distributions of sample maxima, to compute the extremal index, to describe the behavior of extremal clusters, and to derive central limit theorems, as we will do in this section. We will also apply the asymptotic point process results to calculate the limit distributions of the normalized sample autocovariance and autocorrelation functions of the genOU and the increments of the IgenOU process in Section 4.3. The theory that we use goes back to Davis and Hsing (1995) and Davis and Mikosch (1998).

We continue with the definition of a point process. Let the state space $\mathcal{S}$ be $[0,\infty) \times \overline{\mathbb{R}} \setminus \{0\}$, where $\overline{\mathbb{R}} = \mathbb{R} \cup \{-\infty\} \cup \{+\infty\}$. Furthermore, $M_P(\mathcal{S})$ is the class of point measures on $\mathcal{S}$, where $M_P(\mathcal{S})$ is equipped with the metric $\rho$ that generates the topology of vague convergence. The space $(M_P(\mathcal{S}), \rho)$ is a complete and separable metric space with Borel $\sigma$-field $\mathcal{M}_P(\mathcal{S})$. A *point process* in $\mathcal{S}$ is a measurable map from a probability space $(\Omega, \mathcal{A}, \mathbb{P})$ into $(M_P(\mathcal{S}), \mathcal{M}_P(\mathcal{S}))$. A typical example of a point process is a *Poisson random measure*, that is, given a Radon measure $\vartheta$ on $\mathcal{B}(\mathcal{S})$, a point process $\kappa$ is called Poisson random measure with intensity measure $\vartheta$, denoted by $\text{PRM}(\vartheta)$, if

(a) $\kappa(A)$ is Poisson distributed with mean $\vartheta(A)$ for every $A \in \mathcal{B}(\mathcal{S})$,
(b) for all mutually disjoint sets $A_1, \ldots, A_n \in \mathcal{B}(\mathcal{S})$, $\kappa(A_1), \ldots, \kappa(A_n)$ are independent.

More about point processes can be found in Daley and Vere-Jones (2003) and Kallenberg (1997). In our setup we obtain the following result:

**Theorem 4.1 (Point process convergence).** *Let $(V_t)_{t\geq 0}$ be a genOU process satisfying (A). Further, let $(H_k)_{k\in\mathbb{N}}$ and $(I_k)_{k\in\mathbb{N}}$, respectively, be the stationary processes in (1.3) and (1.4). Let $0 < a_n \uparrow \infty$ be a sequence of constants such that*

$$\lim_{n\to\infty} n\mathbb{P}(V_0 > a_n x) = x^{-\alpha} \qquad \text{for } x > 0.$$

(a) *Suppose (B) is satisfied and $\mathcal{A}(a_n)$ holds for $(H_k)_{k\in\mathbb{N}}$. Then*

$$\sum_{k=1}^{\infty} \varepsilon_{(k/n, a_n^{-1} H_k)} \overset{n\to\infty}{\Longrightarrow} \sum_{k=1}^{\infty} \sum_{j=0}^{\infty} \varepsilon_{(s_k^{(1)}, Q_{kj}^{(1)} P_k^{(1)})} \qquad \text{in } \mathcal{M}_P(\mathcal{S}),$$

*where $\sum_{k=1}^{\infty} \varepsilon_{(s_k^{(1)}, P_k^{(1)})}$ is $\mathrm{PRM}(\vartheta)$ with*

$$\vartheta(\mathrm{d}t \times \mathrm{d}x) = \mathrm{d}t \times \alpha h \mathbb{E}\Big(\sup_{0\leq s\leq 1} e^{-\alpha\xi_s} - \sup_{s\geq 1} e^{-\alpha\xi_s}\Big)^{+} x^{-\alpha-1} \mathbf{1}_{(0,\infty)}(x) \, \mathrm{d}x.$$

*Moreover, $\sum_{j=0}^{\infty} \varepsilon_{Q_{kj}^{(1)}}$ for $k \in \mathbb{N}$ are i.i.d. point processes independent of $\sum_{k=1}^{\infty} \varepsilon_{(s_k^{(1)}, P_k^{(1)})}$ with $0 \leq Q_{kj}^{(1)} \leq 1$, and for each $k$ exactly one $Q_{kj}^{(1)}$ is equal to 1, and $\mathbb{P}(Q_{kj}^{(1)} = 0) < 1$ for $j \in \mathbb{N}$. The sequence $(Q_{kj}^{(1)})_{j\in\mathbb{N}_0}$ is a.s. unique.*

(b) *Suppose (C) is satisfied and $\mathcal{A}(a_n^{1/2})$ holds for $(I_k)_{k\in\mathbb{N}}$. Then*

$$\sum_{k=1}^{\infty} \varepsilon_{(k/n, a_n^{-1/2} I_k)} \overset{n\to\infty}{\Longrightarrow} \sum_{k=1}^{\infty} \sum_{j=0}^{\infty} \varepsilon_{(s_k^{(2)}, Q_{kj}^{(2)} P_k^{(2)})} \qquad \text{in } \mathcal{M}_P(\mathcal{S}),$$

*where $\sum_{k=1}^{\infty} \varepsilon_{(s_k^{(2)}, P_k^{(2)})}$ is $\mathrm{PRM}(\vartheta)$ with*

$$\vartheta(\mathrm{d}t \times \mathrm{d}x)$$
$$= \mathrm{d}t \times 2\alpha \mathbb{E}\bigg(\bigg[\Big(\int_0^h e^{-\xi_{t-}/2}\,\mathrm{d}L_t\Big)^{+}\bigg]^{2\alpha}$$
$$\qquad - \max_{k\geq 2}\bigg[\Big(\int_{(k-1)h}^{kh} e^{-\xi_{t-}/2}\,\mathrm{d}L_t\Big)^{+}\bigg]^{2\alpha}\bigg)^{+} x^{-2\alpha-1} \mathbf{1}_{(0,\infty)}(x)\,\mathrm{d}x.$$

*Furthermore, $\sum_{j=0}^{\infty} \varepsilon_{Q_{kj}^{(2)}}$ for $k \in \mathbb{N}$ are i.i.d. point processes independent of $\sum_{k=1}^{\infty} \varepsilon_{(s_k^{(2)}, P_k^{(2)})}$ with $|Q_{kj}^{(2)}| \leq 1$, and for each $k$ exactly one $Q_{kj}^{(2)}$ is equal to 1, and $\mathbb{P}(Q_{kj}^{(2)} = 0) < 1$ for $j \in \mathbb{N}$. The sequence $(Q_{kj}^{(2)})_{j\in\mathbb{N}_0}$ is a.s. unique.*

### 4.1. Extremal behavior

We obtain from Theorem 4.1 the limit behavior of the sequence of partial maxima of the continuous-time process $(V_t)_{t\geq 0}$.

**Proposition 4.2.** *Let the assumptions of Theorem 4.1(a) hold. Define $M(n) := \sup_{0 \leq t \leq n} V_t$ for $n > 0$. Then*

$$\lim_{n \to \infty} \mathbb{P}(a_n^{-1} M(n) \leq x) = \exp\left(-\mathbb{E}\left(\sup_{0 \leq s \leq 1} e^{-\alpha \xi_s} - \sup_{s \geq 1} e^{-\alpha \xi_s}\right)^+ x^{-\alpha}\right) \quad \text{for } x > 0.$$

***Definition 4.3.*** *Let $(X_t)_{t \geq 0}$ be a stationary process. Define for $h > 0$ the sequence $M_k(h) = \sup_{(k-1)h \leq t \leq kh} X_t$, $k \in \mathbb{N}$. If there exist sequences of constants $a_n^{(h)} > 0, b_n^{(h)} \in \mathbb{R}$, a constant $\theta(h) \in [0,1]$ and a non-degenerate distribution function $G$ such that*

$$\lim_{n \to \infty} n\mathbb{P}(M_1(h) > a_n^{(h)} x + b_n^{(h)}) = -\log(G(x)) \quad \text{and}$$

$$\lim_{n \to \infty} \mathbb{P}\left(\max_{k=1,\ldots,n} M_k(h) \leq a_n^{(h)} x + b_n^{(h)}\right) = G(x)^{\theta(h)} \qquad \forall x \text{ in the support of } G,$$

*then we call the function $\theta : (0, \infty) \to [0, 1]$ an extremal index function.*

For fixed $h$ the constant $\theta(h)$ is the *extremal index* of $(M_k(h))_{k \in \mathbb{N}}$ (see Leadbetter (1983), page 67) which is a measure of extremal clusters. The reciprocal of the extremal index can be interpreted as the mean of the cluster size of high-level exceedances: the value 1 reflects no clusters, and values less than 1 reflect clusters.

**Corollary 4.4.**

(a) *Let the assumptions of Theorem 4.1(a) hold. Then*

$$\theta(h) = h \frac{\mathbb{E}(\sup_{0 \leq s \leq 1} e^{-\alpha \xi_s} - \sup_{s \geq 1} e^{-\alpha \xi_s})^+}{\mathbb{E}(\sup_{0 \leq s \leq h} e^{-\alpha \xi_s})} \qquad \text{for } h > 0$$

*is the extremal index function of $(V_t)_{t \geq 0}$.*

(b) *Let the assumptions of Theorem 4.1(b) hold. Then*

$$\theta(h) = \frac{\mathbb{E}([(\int_0^h e^{-\xi_{t-}/2} \, dL_t)^+]^{2\alpha} - \max_{k \geq 2}[(\int_{(k-1)h}^{kh} e^{-\xi_{t-}/2} \, dL_t)^+]^{2\alpha})^+}{\mathbb{E}([(\int_0^h e^{-\xi_{t-}/2} \, dL_t)^+]^{2\alpha})}$$

*is the extremal index of $(I_k)_{k \in \mathbb{N}}$.*

One conclusion is that the processes $(V_t)_{t \geq 0}$, $(H_k)_{k \in \mathbb{N}}$ and $(I_k)_{k \in \mathbb{N}}$ exhibit extremal clusters.

In Fasen *et al.* (2006), the extremal behavior of a COGARCH(1,1) process driven by a compound Poisson process was derived. The next lemma shows that their Theorem 4.5, which says

$$\lim_{n \to \infty} \mathbb{P}(a_n^{-1} M(n) \leq x) = \exp\left(\mu(\mathbb{E}(e^{-\alpha c \Gamma_1}))^{-1} \mathbb{E}\left(1 - \sup_{s \geq \Gamma_1} e^{-\alpha \xi_s}\right)^+ x^{-\alpha}\right) \quad \text{for } x > 0$$

with the notation of Lemma 4.5 below, and our Proposition 4.2 are consistent.

**Lemma 4.5.** *Let $(V_t)_{t\geq 0}$ be the volatility process of the $\mathrm{COGARCH}(1,1)$ model in (2.7) satisfying (2.8). Let $(L_t)_{t\geq 0}$ be a compound Poisson process with jump arrivals $(\Gamma_k)_{k\in\mathbb{N}}$ and intensity $\mu$. Then*

$$\mathbb{E}\Big(\sup_{0\leq s\leq 1} e^{-\alpha\xi_s} - \sup_{s\geq 1} e^{-\alpha\xi_s}\Big)^+ = \mu(\mathbb{E}(e^{-\alpha c\Gamma_1}))^{-1}\mathbb{E}\Big(1 - \sup_{s\geq \Gamma_1} e^{-\alpha\xi_s}\Big)^+. \qquad (4.1)$$

## 4.2. Asymptotic behavior of the IgenOU process

The last conclusion of Theorem 4.1 is a central limit theorem for $(I_t^*)_{t\geq 0}$.

**Proposition 4.6.** *Let $(I_t^*)_{t\geq 0}$ be the IgenOU process in (1.2), and $(I_k)_{k\in\mathbb{N}}$ as in (1.4) satisfies (A) and (C). Let $0 < a_n \uparrow \infty$ be a sequence of constants such that*

$$\lim_{n\to\infty} n\mathbb{P}(V_0 > a_n x) = x^{-\alpha} \qquad \text{for } x > 0.$$

(a) *If $\alpha \in (0, 0.5)$ and $(I_k)_{k\in\mathbb{N}}$ satisfies $\mathcal{A}(a_n^{1/2})$, then*

$$t^{-1/(2\alpha)} I_t^* \stackrel{t\to\infty}{\Longrightarrow} S, \qquad \text{where } S \text{ is } (2\alpha)\text{-stable.}$$

(b) *If $\alpha \in (0.5, 1)$, $(I_k)_{k\in\mathbb{N}}$ satisfies $\mathcal{A}(a_n^{1/2})$ and*

$$\lim_{\epsilon\downarrow 0}\limsup_{n\to\infty} \mathrm{Var}\left(n^{-1/(2\alpha)} \sum_{k=1}^n I_k \mathbf{1}_{\{|I_k|\leq \epsilon n^{1/(2\alpha)}\}}\right) = 0, \qquad (4.2)$$

*then*

$$t^{-1/(2\alpha)}(I_t^* - \mathbb{E}(I_t^*)) \stackrel{t\to\infty}{\Longrightarrow} S, \qquad \text{where } S \text{ is } (2\alpha)\text{-stable.}$$

(c) *If $\alpha > 1$ and $(I_k)_{k\in\mathbb{N}}$ is exponentially $\alpha$-mixing, then*

$$t^{-1/2}(I_t^* - \mathbb{E}(I_t^*)) \stackrel{t\to\infty}{\Longrightarrow} \mathcal{N},$$

*where $\mathcal{N}$ is normal distributed with $\mathbb{E}(\mathcal{N}) = 0$ and $\mathrm{Var}(\mathcal{N}) = \mathrm{Var}(I_1^*)$.*

This proposition is a consequence of Davis and Hsing (1995), Theorem 3.1, and Ibragimov and Linnik (1971), Theorem 18.5.3.

**Remark 4.7.**

(i) Condition (4.2) is satisfied if $(\xi_t, \eta_t, L_t)_{t\geq 0} \stackrel{d}{=} (\xi_t, \eta_t, -L_t)_{t\geq 0}$, since $I_k$ then is symmetric. Thus, (4.2) stems from the uncorrelation of $(I_k)_{k\in\mathbb{N}}$ and Karamata's theorem (see Feller (1971), VIII.9, Theorem 1). A necessary but insufficient condition of $(\xi_t, \eta_t, L_t)_{t\geq 0} \stackrel{d}{=} (\xi_t, \eta_t, -L_t)_{t\geq 0}$ is $L$ symmetric. For example, let $L_1$ be symmetric and independent of the subordinator $\eta$. Then $(L_t, \eta_t, L_t)_{t\geq 0} \stackrel{d}{=} (-L_t, \eta_t, -L_t)_{t\geq 0}$, but the distribution differs from the distribution of $(L_t, \eta_t, -L_t)_{t\geq 0}$.

(ii) Let us consider the COGARCH$(1,1)$ model in Example 2.4 and suppose that $L_1$ has a symmetric distribution. Then $(\xi_t, t, L_t)_{t\geq 0} \stackrel{d}{=} (\xi_t, t, -L_t)_{t\geq 0}$, and (4.2) holds. In particular (4.2) also holds for the Nelson diffusion model.

(iii) The boundary cases $\alpha = 0.5, 1$ are here neglected, since the analysis is tedious and lengthy, and it does not lead to interesting statistical insight.

### 4.3. Convergence of the sample autocovariances

The next section is devoted to the asymptotic behavior of the sample autocovariance and autocorrelation function of $(V_t)_{t\geq 0}$, and $(I_k)_{k\in\mathbb{N}}$.

**Theorem 4.8.** *Let $(V_t)_{t\geq 0}$ be a genOU process satisfying* (A) *and* (B). *Suppose $(V_t)_{t\geq 0}$ is exponentially $\alpha$-mixing. Further, let $\gamma_V(t) = \mathbb{E}(V_0 V_t)$ and $\rho_V(t) = \gamma_V(t)/\gamma_V(0)$ for $t > 0$. Define for $h > 0$ the empirical versions*

$$\gamma_{n,V}(lh) = \frac{1}{n}\sum_{k=1}^{n-l} V_{kh} V_{(k+l)h} \quad and \quad \rho_{n,V}(lh) = \gamma_{n,V}(lh)/\gamma_{n,V}(0) \qquad for\ l \in \mathbb{N}_0.$$

(a) *If $\alpha \in (0, 2)$, then*

$$(n^{1-2/\alpha}\gamma_{n,V}(lh))_{l=0,\ldots,m} \stackrel{n\to\infty}{\Longrightarrow} (S_l^{(1)})_{l=0,\ldots,m}, \tag{4.3}$$

$$(\rho_{n,V}(lh))_{l=1,\ldots,m} \stackrel{n\to\infty}{\Longrightarrow} (S_l^{(1)}/S_0^{(1)})_{l=1,\ldots,m}, \tag{4.4}$$

*where the vector $(S_0^{(1)}, \ldots, S_m^{(1)})$ is jointly $(\alpha/2)$-stable in $\mathbb{R}^{m+1}$.*

(b) *If $\alpha \in (2, 4)$ and $d > 4$ in condition* (B), *then*

$$(n^{1-2/\alpha}(\gamma_{n,V}(lh) - \gamma_V(lh)))_{l=0,\ldots,m} \stackrel{n\to\infty}{\Longrightarrow} (S_l^{(2)})_{l=0,\ldots,m}, \tag{4.5}$$

$$(n^{1-2/\alpha}(\rho_{n,V}(lh) - \rho_V(lh)))_{l=1,\ldots,m} \stackrel{n\to\infty}{\Longrightarrow} \gamma_V^{-1}(0)(S_l^{(2)} - \rho_V(lh)S_0^{(2)})_{l=1,\ldots,m}, \tag{4.6}$$

*where $(S_0^{(2)}, \ldots, S_m^{(2)})$ is jointly $(\alpha/2)$-stable in $\mathbb{R}^{m+1}$.*

(c) *If $\alpha > 4$, then (4.5) and (4.6) hold with normalization $n^{1/2}$, where the limit $(S_1^{(3)}, \ldots, S_m^{(3)})$ is multivariate normal with mean zero, covariance matrix*

$$\left(\sum_{k=-\infty}^{\infty} \mathrm{Cov}(V_0 V_{ih}, V_{kh} V_{(k+j)h})\right)_{i,j=1,\ldots,m} \quad and \quad S_0^{(3)} = \mathbb{E}(V_0^2).$$

*Remark 4.9.*

(i) The stable random vector $(S_0^{(1)}, \ldots, S_m^{(1)})$ is a functional of the limit point process based on $(V_{kh})_{k\in\mathbb{N}}$ in (B.12). The explicit representation of $(S_0^{(1)}, \ldots, S_m^{(1)})$ is given in (B.19). Similarly we can derive the representation of $(S_0^{(2)}, \ldots, S_m^{(2)})$.

(ii) If $\alpha \in (0,2)$, the autocovariance function does not exist. Hence, $\gamma_{n,V}$ and $\rho_{n,V}$ are not consistent estimators.

(iii) For $\alpha > 2$ the sample autocovariance function is a consistent estimator, where for $\alpha \in (2,4)$ the convergence rate $n^{1-2/\alpha}$ will be faster, if $\alpha$ increases. The convergence to an infinite variance stable distribution in (b) and the slower convergence rate than in (c) cause the confidence bands in (b) to be wider than in (c).

(iv) The mean corrected versions of the sample and the autocovariance function can also be considered; the limit theory does not change.

(v) The proof of Theorem 4.8 shows that (b) is valid under more general assumptions. Let $(\widetilde{V}_k)_{k \in \mathbb{N}}$ be the stationary solution of the stochastic recurrence equation $\widetilde{V}_{k+1} = \widetilde{A}_k \widetilde{V}_k + \widetilde{B}_k$, where $(\widetilde{A}_k, \widetilde{B}_k)_{k \in \mathbb{N}}$ is an i.i.d. sequence, and $(\widetilde{A}_k, \widetilde{B}_k)$ is independent of $\widetilde{V}_k$. Let $(\widetilde{A}_k \widetilde{B}_k \widetilde{V}_k)_{k \in \mathbb{N}}$ be exponentially $\alpha$-mixing. Furthermore, we suppose that the finite dimensional distributions of $(\widetilde{V}_k)_{k \in \mathbb{N}}$ are multivariate regularly varying of index $\alpha \in (2,4)$, and $\mathbb{E}|\widetilde{A}_k|^d < \infty$ and $\mathbb{E}|\widetilde{B}_k|^d < \infty$ for some $d > 4$. If, finally, $(\widetilde{V}_k)_{k \in \mathbb{N}}$ satisfies (B.12) with $V_k$ replaced by $\widetilde{V}_k$, then Theorem 4.8(b) holds.

The following result is a straightforward conclusion of Theorem 3.3, Theorem 4.1 and Davis and Mikosch (1998), Theorem 3.5.

**Theorem 4.10.** *Let $(I_k)_{k \in \mathbb{N}}$ be the stationary process in (1.4) satisfying (A) and (C). Suppose $(I_k)_{k \in \mathbb{N}}$ is exponentially $\alpha$-mixing. Further, let $\gamma_I(l) = \mathbb{E}(I_1 I_{1+l})$ and $\rho_I(l) = \gamma_I(l)/\gamma_I(0)$ for $l \in \mathbb{N}$. Define the empirical versions*

$$\gamma_{n,I}(l) = \frac{1}{n} \sum_{k=1}^{n-l} I_k I_{k+l} \quad and \quad \rho_{n,I}(l) = \gamma_{n,I}(l)/\gamma_{n,I}(0) \qquad for\ l \in \mathbb{N}_0.$$

(a) *If $\alpha \in (0,1)$ then*

$$(n^{1-1/\alpha} \gamma_{n,I}(l))_{l=0,\ldots,m} \stackrel{n \to \infty}{\Longrightarrow} (S_l^{(1)})_{l=0,\ldots,m}, \tag{4.7}$$

$$(\rho_{n,I}(l))_{l=1,\ldots,m} \stackrel{n \to \infty}{\Longrightarrow} (S_l^{(1)}/S_0^{(1)})_{l=1,\ldots,m}, \tag{4.8}$$

*where the vector $(S_0^{(1)}, \ldots, S_m^{(1)})$ is jointly $\alpha$-stable in $\mathbb{R}^{m+1}$.*

(b) *If $\alpha \in (1,2)$ and*

$$\lim_{\epsilon \downarrow 0} \limsup_{n \to \infty} \operatorname{Var} \left( n^{-1/\alpha} \sum_{i=1}^{n-l} I_i I_{i+l} \mathbf{1}_{\{|I_i I_{i+l}| \leq n^{1/\alpha} \epsilon\}} \right) = 0, \qquad l = 0, \ldots, m, \tag{4.9}$$

*then*

$$(n^{1-1/\alpha}(\gamma_{n,I}(l) - \gamma_I(l)))_{l=0,\ldots,m} \stackrel{n \to \infty}{\Longrightarrow} (S_l^{(2)})_{l=0,\ldots,m}, \tag{4.10}$$

$$(n^{1-1/\alpha}(\rho_{n,I}(l) - \rho_I(l)))_{l=1,\ldots,m} \stackrel{n \to \infty}{\Longrightarrow} \gamma_I^{-1}(0)(S_l^{(2)} - \rho_I(l)S_0^{(2)})_{l=1,\ldots,m}, \tag{4.11}$$

where $(S_0^{(2)}, \ldots, S_m^{(2)})$ is jointly $\alpha$-stable in $\mathbb{R}^{m+1}$.

(c) If $\alpha > 2$ then (4.10) and (4.11) hold with normalization $n^{1/2}$, where the limit $(S_1^{(3)}, \ldots, S_m^{(3)})$ is multivariate normal with mean zero, covariance matrix

$$\left(\sum_{k=-\infty}^{\infty} \mathrm{Cov}(I_1 I_{1+i}, I_k I_{k+j})\right)_{i,j=1,\ldots,m} \quad \text{and} \quad S_0^{(3)} = \mathbb{E}(I_1^2).$$

As in Remark 4.7, a sufficient condition for (4.9) is $(\xi_t, \eta_t, L_t)_{t \geq 0} \stackrel{d}{=} (\xi_t, \eta_t, -L_t)_{t \geq 0}$.

## Appendix A: Proofs of Sections 2 and 3

**Remark A.1.** Let condition (B) be satisfied.

(i) By Sato (1999), Lemma 26.4, we know that $\Psi_\xi$ is strictly convex and continuous. Hence, $\Psi_\xi(\alpha) = 0$ implies that there exists a $0 < \widetilde{\alpha} < \alpha$ such that

$$\Psi_\xi(\widetilde{\alpha}) < 0. \tag{A.1}$$

(ii) The process $(\mathrm{e}^{-v\xi_t - t\Psi_\xi(v)})_{t \geq 0}$ is a martingale for every $v \in \mathbb{R}$ where $|\Psi_\xi(v)| < \infty$. A conclusion of Doob's martingale inequality (cf. Revuz and Yor (2001), page 54) is that $K_{v,h} := \mathbb{E}(\sup_{0 \leq t \leq h} \mathrm{e}^{-v\xi_t}) < \infty$ for $h > 0$, and hence, by the independent and stationary increments of a Lévy process

$$\mathbb{E}\left(\sup_{(k-1)h \leq t \leq kh} \mathrm{e}^{-v\xi_t}\right) = \mathbb{E}\left(\sup_{(k-1)h \leq t \leq kh} \mathrm{e}^{-v(\xi_t - \xi_{(k-1)h})}\right) \mathbb{E}(\mathrm{e}^{-v\xi_{(k-1)h}}) \\ \leq K_{v,h} \mathrm{e}^{\Psi_\xi(v)(k-1)h} < \infty. \tag{A.2}$$

(iii) Let $0 < u \leq d$. Then there exists a $0 < v < \min(u, \alpha)$ such that $((\mathrm{e}^{-\xi_t} \int_0^t \mathrm{e}^{\xi_{s-}} \, \mathrm{d}\eta_s)^v)_{t \geq 0}$ is a positive submartingale. Doob's submartingale inequality, Hölder's inequality and (2.2) result in

$$\mathbb{E}\left(\sup_{0 \leq t \leq h} \mathrm{e}^{-\xi_t} \int_0^t \mathrm{e}^{\xi_{s-}} \, \mathrm{d}\eta_s\right)^u \leq \widetilde{K}_{u,h} \left(\mathbb{E}\left(\mathrm{e}^{-\xi_h} \int_0^h \mathrm{e}^{\xi_{s-}} \, \mathrm{d}\eta_s\right)^d\right)^{u/d} < \infty \tag{A.3}$$

for some constant $\widetilde{K}_{u,h} > 0$.

**Proof of Proposition 2.1.**

(a) Condition (A) follows by Lindner and Maller (2005), Theorem 4.5. Hence, it remains only to prove (2.2). Since $\eta$ is a subordinator we have

$$\mathbb{E}\left|\mathrm{e}^{-\xi_h} \int_0^h \mathrm{e}^{\xi_{s-}} \, \mathrm{d}\eta_s\right|^d \leq \mathbb{E}\left|\sup_{0 \leq s \leq h} \mathrm{e}^{-(\xi_h - \xi_s)} \eta_h\right|^d.$$

Applying Hölder's inequality and (A.2) we obtain

$$\mathbb{E}\left|e^{-\xi_h}\int_0^h e^{\xi_{s-}}\,d\eta_s\right|^d \leq C(\mathbb{E}(e^{-pd\xi_h}))^{1/p}(\mathbb{E}|\eta_h|^{qd})^{1/q} < \infty. \tag{A.4}$$

(b) By (a) we have only to check (2.3). We assume $\mathbb{E}(L_1) = 0$, or else we decompose $L$ into two independent Lévy processes where one process has mean 0 and the other is a drift term. Then $(\int_0^u e^{-\xi_{t-}/2}(\int_0^{t-} e^{\xi_{s-}}\,d\eta_s)^{1/2}\,dL_t)_{u\geq 0}$ is a local martingale by Protter (2004), Theorem 29, page 173. Further, we define $\widetilde{d} = \max\{1,d\}$. By the Burkholder–Gundy inequality (cf. Liptser and Shiryayev (1989), page 75), Hölder's inequality and (A.3) we obtain

$$\mathbb{E}\left|\int_0^h e^{-\xi_{t-}/2}\left(\int_0^{t-} e^{\xi_{s-}}\,d\eta_s\right)^{1/2}dL_t\right|^{2\widetilde{d}}$$

$$\leq K_1 \mathbb{E}\left|\int_0^h e^{-\xi_{t-}}\left(\int_0^{t-} e^{\xi_{s-}}\,d\eta_s\right)d[L,L]_t\right|^{\widetilde{d}}$$

$$\leq K_2\left(\mathbb{E}\left|\sup_{0\leq t\leq h} e^{-\xi_{t-}}\int_0^{t-} e^{\xi_{s-}}\,d\eta_s\right|^{\widetilde{d}p_2}\right)^{1/p_2}(\mathbb{E}|L_h|^{2\widetilde{d}q_2})^{1/q_2} < \infty.$$

The finiteness of the first factor is again a conclusion of (A.3), (A.4) and (2.5).

Similarly, we can prove that $\mathbb{E}|\int_0^h e^{-\xi_{t-}/2}\,dL_t|^{2\widetilde{d}} < \infty$. □

**Proof of Theorem 3.3.**

(a) *Step 1.* Let $\overline{\mathbf{V}} = (e^{-\xi_t}V_0)_{0\leq t\leq 1}$. Since the tail of the probability distribution of $V_0$ is regularly varying, the process $\overline{\mathbf{V}}$ is a regularly varying process by Hult and Lindskog (2007), Theorem 3.1, and

$$\frac{\mathbb{P}(|\overline{\mathbf{V}}|_\infty > ux, \overline{\mathbf{V}}/|\overline{\mathbf{V}}|_\infty \in \cdot)}{\mathbb{P}(|\overline{\mathbf{V}}|_\infty > u)} \stackrel{u\to\infty}{\Longrightarrow} x^{-\alpha}\frac{\mathbb{E}(|\mathbf{U}|_\infty^\alpha \mathbf{1}\{\mathbf{U}/|\mathbf{U}|_\infty \in \cdot\})}{\mathbb{E}|\mathbf{U}|_\infty^\alpha} \qquad \text{on } \mathcal{B}(\mathbb{S}_\mathbb{D}). \tag{A.5}$$

*Step 2.* We will show that $\mathbf{V}$ is a regularly varying process. By (A.3) we know that $\mathbb{E}|\mathbf{V}-\overline{\mathbf{V}}|_\infty^d < \infty$. Markov's inequality, regular variation of $\mathbb{P}(|\overline{\mathbf{V}}|_\infty > u)$ as $u\to\infty$ (by Step 1) and Potter's theorem (cf. Bingham *et al.* (1987), Theorem 1.5.6) give

$$\frac{\mathbb{P}(|\mathbf{V}-\overline{\mathbf{V}}|_\infty > u)}{\mathbb{P}(|\overline{\mathbf{V}}|_\infty > u)} \leq \frac{u^{-d}}{\mathbb{P}(|\overline{\mathbf{V}}|_\infty > u)}\mathbb{E}|\mathbf{V}-\overline{\mathbf{V}}|_\infty^d \longrightarrow 0 \qquad \text{as } u\to\infty.$$

Hence, as in Jessen and Mikosch (2006), Lemma 3.12,

$$\frac{\mathbb{P}(|\mathbf{V}|_\infty > ux, \mathbf{V}/|\mathbf{V}|_\infty \in \cdot)}{\mathbb{P}(|\mathbf{V}|_\infty > u)} \sim \frac{\mathbb{P}(|\overline{\mathbf{V}}|_\infty > ux, \overline{\mathbf{V}}/|\overline{\mathbf{V}}|_\infty \in \cdot)}{\mathbb{P}(|\overline{\mathbf{V}}|_\infty > u)} \qquad \text{as } u\to\infty.$$

With Step 1, part (a) follows.

(b) Analogous to (a), we have $\mathbf{V}^{(kh)} = (V_t)_{0 \leq t \leq kh}$ is a regularly varying process in $\mathbb{D}[0, kh]$. The functional $T: \mathbb{D}[0, kh] \setminus \{\mathbf{0}\} \to \mathbb{R}^k$ with

$$\mathbf{x} = (x_t)_{0 \leq t \leq kh} \mapsto \left( \sup_{0 \leq s \leq h} |x_s|, \ldots, \sup_{(k-1)h \leq s \leq kh} |x_s| \right)$$

is continuous, $|\mathbf{x}|_\infty = |T(\mathbf{x})|_\infty$ and $T(\lambda \mathbf{x}) = \lambda T(\mathbf{x})$ for $\lambda > 0, \mathbf{x} \in \mathbb{D}[0, kh]$. Let $\mathbf{U}^{(kh)} = (\mathrm{e}^{-\xi_t})_{0 \leq t \leq kh}$. Then

$$\mathbf{H}_k = T(\mathbf{V}^{(kh)}), \qquad |\mathbf{H}_k|_\infty = |\mathbf{V}^{(kh)}|_\infty \quad \text{and} \quad T(\mathbf{V}^{(kh)}/|\mathbf{V}^{(kh)}|_\infty) = \mathbf{H}_k/|\mathbf{H}_k|_\infty,$$

and similarly

$$\mathbf{m}_k = T(\mathbf{U}^{(kh)}), \quad |\mathbf{m}_k|_\infty = |\mathbf{U}^{(kh)}|_\infty \quad \text{and} \quad T(\mathbf{U}^{(kh)}/|\mathbf{U}^{(kh)}|_\infty) = \mathbf{m}_k/|\mathbf{m}_k|_\infty.$$

We conclude by the continuous mapping theorem (cf. Billingsley (1999), Theorem 2.7, and Hult and Lindskog (2005), Theorem 8, for regularly varying stochastic processes) and (a) that on $\mathcal{B}(\mathbb{S}_{k-1})$,

$$\frac{\mathbb{P}(|\mathbf{H}_k|_\infty > ux, \mathbf{H}_k/|\mathbf{H}_k|_\infty \in \cdot)}{\mathbb{P}(|\mathbf{H}_k|_\infty > u)} = \frac{\mathbb{P}(|\mathbf{V}^{(kh)}|_\infty > ux, T(\mathbf{V}^{(kh)}/|\mathbf{V}^{(kh)}|_\infty) \in \cdot)}{\mathbb{P}(|\mathbf{V}^{(kh)}|_\infty > u)}$$

$$\stackrel{u \to \infty}{\Longrightarrow} \frac{\mathbb{E}(|\mathbf{U}^{(kh)}|_\infty^\alpha \mathbf{1}\{T(\mathbf{U}^{(kh)}/|\mathbf{U}^{(kh)}|_\infty) \in \cdot\})}{\mathbb{E}|\mathbf{U}^{(kh)}|_\infty^\alpha},$$

which gives the desired result.

(c) We define

$$I_k^{(1)} = \sqrt{V_0} \int_{(k-1)h}^{kh} \mathrm{e}^{-\xi_{t-}/2} \, \mathrm{d}L_t \quad \text{and} \quad I_k^{(2)} = I_k - I_k^{(1)} = \int_{(k-1)h}^{kh} R_t \, \mathrm{d}L_t,$$

where $R_t = \sqrt{V_{t-}} - \mathrm{e}^{-\xi_{t-}/2}\sqrt{V_0}$. Note that

$$R_t^2 \leq \mathrm{e}^{-\xi_{t-}} \int_0^{t-} \mathrm{e}^{\xi_{s-}} \, \mathrm{d}\eta_s.$$

We assume $\mathbb{E}(L_1) = 0$, or else we decompose $L$ into two independent Lévy processes where one process has mean 0 and the other is a drift term. Then $(\int_0^u R_t \, \mathrm{d}L_t)_{u \geq 0}$ is a local martingale by Protter (2004), Theorem 29, page 173. Further, we define $\widetilde{d} = \max\{1, d\}$. By the Burkholder–Gundy inequality (cf. Liptser and Shiryayev (1989), page 75) and (A.3) we obtain

$$\mathbb{E}|I_k^{(2)}|^{2\widetilde{d}} \leq K_1 \mathbb{E}\left( \int_{(k-1)h}^{kh} R_t^2 \, \mathrm{d}[L, L]_t \right)^{\widetilde{d}}$$

$$\leq K_1 \mathbb{E}\bigg(\int_{(k-1)h}^{kh} \bigg(\mathrm{e}^{-\xi_{t-}} \int_0^{t-} \mathrm{e}^{\xi_{s-}}\,\mathrm{d}\eta_s\bigg)\mathrm{d}[L,L]_t\bigg)^{\widetilde{d}} \tag{A.6}$$

$$\leq K_2 \mathbb{E}\bigg|\int_{(k-1)h}^{kh} \mathrm{e}^{-\xi_{t-}/2}\bigg(\int_0^{t-} \mathrm{e}^{\xi_{s-}}\,\mathrm{d}\eta_s\bigg)^{1/2}\mathrm{d}L_t\bigg|^{2\widetilde{d}} < \infty,$$

where the finiteness follows from (2.3). Thus, the classical result of Breiman (1965), and Klüppelberg *et al.* (2006), Lemma 2, leads to

$$\mathbb{P}(I_1 > x) \sim \mathbb{P}(I_1^{(1)} > x) \sim \mathbb{E}\bigg[\bigg(\int_0^h \mathrm{e}^{-\xi_{t-}/2}\,\mathrm{d}L_t\bigg)^+\bigg]^{2\alpha} \mathbb{P}(\sqrt{V_0} > x) \qquad \text{as } x \to \infty.$$

With Basrak *et al.* (2002), Proposition A.1, which is a multivariate version of Breiman's result, we can extend this result to the multivariate case of $\mathbf{I}_k$. □

**Proof of Proposition 3.4.**

(a) Follows along the same lines as the proof of Theorem 4.3 in Masuda (2004), where the result was derived for the classical Ornstein–Uhlenbeck process.

(b) follows from (a).

(c) follows from (a) and Genon-Catalot *et al.* (2000), Proposition 3.1. The arguments are the same as in Theorem 3.1 and Proposition 3.2 of Genon-Catalot *et al.* (2000), who investigate a slightly different model. □

## Appendix B: Proofs of Section 4

The proof of Theorem 4.1 uses the next lemma.

**Lemma B.1.** *Let* $(\xi_t)_{t\geq 0}$ *be a Lévy process satisfying* $\mathbb{E}(\mathrm{e}^{-\alpha\xi_1}) = 1$. *Then*

$$\mathbb{E}\bigg(\sup_{0\leq s\leq h} \mathrm{e}^{-\alpha\xi_s} - \sup_{s\geq h} \mathrm{e}^{-\alpha\xi_s}\bigg)^+$$

$$= h\mathbb{E}\bigg(\sup_{0\leq s\leq 1} \mathrm{e}^{-\alpha\xi_s} - \sup_{s\geq 1} \mathrm{e}^{-\alpha\xi_s}\bigg)^+ \qquad \text{for any } h > 0.$$

**Proof.** For $f:[0,\infty) \to \mathbb{R}$ holds

$$\bigg(\sup_{s\leq kh} f(s) - \sup_{s\geq kh} f(s)\bigg)^+$$

$$= \bigg(\sup_{s\leq h} f(s) - \sup_{s\geq h} f(s)\bigg)^+ + \sum_{l=1}^{k-1}\bigg(\sup_{lh\leq s\leq (l+1)h} f(s) - \sup_{s\geq (l+1)h} f(s)\bigg)^+. \tag{B.1}$$

This can be proved by induction, since

$$\Big(\sup_{s\leq kh} f(s) - \sup_{s\geq kh} f(s)\Big)^+$$
$$= \Big(\sup_{s\leq (k-1)h} f(s) - \sup_{s\geq (k-1)h} f(s)\Big)^+ + \Big(\sup_{(k-1)h\leq s\leq kh} f(s) - \sup_{s\geq kh} f(s)\Big)^+.$$

One way to see this is to distinguish the different cases where the maximum and the second largest maximum of $\sup_{s\leq (k-1)h} f(s)$, $\sup_{(k-1)h\leq s\leq kh} f(s)$ and $\sup_{s\geq kh} f(s)$ lie. Taking the independent increments of $(\xi_t)_{t\geq 0}$ into account, we obtain

$$\mathbb{E}\Big(\sup_{lh\leq s\leq (l+1)h} e^{-\alpha\xi_s} - \sup_{s\geq (l+1)h} e^{-\alpha\xi_s}\Big)^+$$
$$= \mathbb{E}\Big(e^{-\alpha\xi_{lh}}\Big(\sup_{lh\leq s\leq (l+1)h} e^{-\alpha(\xi_s-\xi_{lh})} - \sup_{s\geq (l+1)h} e^{-\alpha(\xi_s-\xi_{lh})}\Big)\Big)^+$$
$$= \mathbb{E}(e^{-\alpha\xi_{lh}})\mathbb{E}\Big(\sup_{lh\leq s\leq (l+1)h} e^{-\alpha(\xi_s-\xi_{lh})} - \sup_{s\geq (l+1)h} e^{-\alpha(\xi_s-\xi_{lh})}\Big)^+.$$

Since $\mathbb{E}(e^{-\alpha\xi_1}) = 1$, by assumption we also have $\mathbb{E}(e^{-\alpha\xi_{lh}}) = 1$. Thus, the stationary increments property of $(\xi_t)_{t\geq 0}$ gives

$$\mathbb{E}\Big(\sup_{lh\leq s\leq (l+1)h} e^{-\alpha\xi_s} - \sup_{s\geq (l+1)h} e^{-\alpha\xi_s}\Big)^+ = \mathbb{E}\Big(\sup_{0\leq s\leq h} e^{-\alpha\xi_s} - \sup_{s\geq h} e^{-\alpha\xi_s}\Big)^+. \quad (B.2)$$

Hence, (B.1) and (B.2) result in

$$\mathbb{E}\Big(\sup_{0\leq s\leq kh} e^{-\alpha\xi_s} - \sup_{s\geq kh} e^{-\alpha\xi_s}\Big)^+ = k\mathbb{E}\Big(\sup_{0\leq s\leq h} e^{-\alpha\xi_s} - \sup_{s\geq h} e^{-\alpha\xi_s}\Big)^+ \quad \text{for } k \in \mathbb{N}.$$

Thus,

$$\mathbb{E}\Big(\sup_{0\leq s\leq hq} e^{-\alpha\xi_s} - \sup_{s\geq hq} e^{-\alpha\xi_s}\Big)^+ = q\mathbb{E}\Big(\sup_{0\leq s\leq h} e^{-\alpha\xi_s} - \sup_{s\geq h} e^{-\alpha\xi_s}\Big)^+ \quad \text{for } q \in \mathbb{Q} \cap \mathbb{R}_+.$$

Since $(\xi_t)_{t\geq 0}$ has a.s. càdlàg paths, the claim follows. □

**Proof of Theorem 4.1.** (a) Davis and Hsing (1995) derived sufficient assumptions for the convergence of point processes formed by a stationary, regularly varying sequence. We apply their Theorem 2.7. This theorem requires that the finite dimensional distributions of $(H_k)_{k\in\mathbb{N}}$ are multivariate regularly varying, which is satisfied by Theorem 3.3(b). Finally, we must check

$$\lim_{l\to\infty}\lim_{n\to\infty} \mathbb{P}\Big(\bigvee_{l\leq k\leq r_n} |H_k| > a_n x \Big| |H_1| > a_n x\Big) = 0 \quad (B.3)$$

for a sequence $r_n = \mathrm{o}(n)$ as $n \to \infty$, and $x > 0$. Thus, define $A_k^* = \sup_{(k-1)h \leq t \leq kh} \mathrm{e}^{-(\xi_t - \xi_h)}$ and $B_k^* = \sup_{(k-1)h \leq t \leq kh} \mathrm{e}^{-\xi_t} \int_h^t \mathrm{e}^{\xi_{s-}} \, \mathrm{d}\eta_s$ for $k \geq 2$. Let $x > 0$ be fixed. Hence, $H_k \leq A_k^* H_1 + B_k^*$ for $k \geq 2$ and $(A_k^*, B_k^*)$ are independent of $H_1$. Then we obtain

$$\mathbb{P}\bigg(\bigvee_{l \leq k \leq r_n} |H_k| > a_n x, |H_1| > a_n x\bigg)$$
$$\leq \sum_{k=l}^{r_n} [\mathbb{P}(B_k^* > a_n x/2) \mathbb{P}(H_1 > a_n x) + \mathbb{P}(A_k^* H_1 > a_n x/2, H_1 > a_n x)] \quad \text{(B.4)}$$
$$\leq r_n \mathbb{P}(H_1 > a_n x/2)^2 + \sum_{k=l}^{r_n} \mathbb{P}(A_k^* H_1 > a_n x/2, H_1 > a_n x).$$

Let $\widetilde{\alpha}$ be given as in (A.1). Markov's inequality and the independence of $H_1$ and $A_k^*$ for $k \geq 2$ lead to

$$\mathbb{P}(A_k^* H_1 \mathbf{1}_{\{H_1 > a_n x\}} > a_n x/2) \leq (a_n x/2)^{-\widetilde{\alpha}} \mathbb{E}(A_k^{*\widetilde{\alpha}}) \mathbb{E}(H_1^{\widetilde{\alpha}} \mathbf{1}_{\{H_1 > a_n x\}}). \quad \text{(B.5)}$$

Since $\mathbb{P}(H_1^{\widetilde{\alpha}} > x)$ is regularly varying with index $\alpha/\widetilde{\alpha} > 1$, we apply Feller (1971), Theorem 1 in Chapter VIII.9, such that

$$\mathbb{E}(H_1^{\widetilde{\alpha}} \mathbf{1}_{\{H_1 > a_n x\}}) \sim \frac{\widetilde{\alpha}}{\alpha - \widetilde{\alpha}} (a_n x)^{\widetilde{\alpha}} \mathbb{P}(H_1 > a_n x) \quad \text{as } n \to \infty. \quad \text{(B.6)}$$

Hence, (B.5), (B.6) and (A.2) result in

$$\mathbb{P}(A_k^* H_1 \mathbf{1}_{\{H_1 > a_n x\}} > a_n x/2) \leq K_1 \mathbb{E}(A_k^{*\widetilde{\alpha}}) \mathbb{P}(H_1 > a_n x)$$
$$\leq K_2 \mathrm{e}^{\Psi_\xi(\widetilde{\alpha})kh} \mathbb{P}(H_1 > a_n x), \quad \text{(B.7)}$$

for $n \geq n_0$ and some constants $K_1, K_2, n_0 > 0$. We conclude from (B.4) and (B.7) that

$$\lim_{l \to \infty} \lim_{n \to \infty} \mathbb{P}\bigg(\bigvee_{l \leq |k| \leq r_n} |H_k| > a_n x \bigg| |H_1| > a_n x\bigg) \leq \lim_{l \to \infty} K_2 \sum_{k=l}^{\infty} (\mathrm{e}^{\Psi_\xi(\widetilde{\alpha})h})^k = 0.$$

Note that (cf. Theorem 2.7 and Lemma 2.9 in Davis and Hsing (1995)) the extremal index $\theta$ of $(H_k)_{k \in \mathbb{N}}$ has value

$$\theta = \lim_{l \to \infty} \lim_{u \to \infty} \frac{\mathbb{P}(|H_1| > u) - \mathbb{P}(\min\{|H_1|, \bigvee_{k=2}^l |H_k|\} > u)}{\mathbb{P}(|H_1| > u)}.$$

By Theorem 3.3(b) the right-hand side is equal to

$$\lim_{l \to \infty} \frac{\mathbb{E}(\sup_{0 \leq s \leq h} \mathrm{e}^{-\alpha \xi_s}) - \mathbb{E}(\min\{\sup_{0 \leq s \leq h} \mathrm{e}^{-\alpha \xi_s}, \sup_{h \leq s \leq lh} \mathrm{e}^{-\alpha \xi_s}\})}{\mathbb{E}(\sup_{0 \leq s \leq h} \mathrm{e}^{-\alpha \xi_s})}$$

$$= \frac{\mathbb{E}(\sup_{0\leq s\leq h} e^{-\alpha\xi_s} - \sup_{s\geq h} e^{-\alpha\xi_s})^+}{\mathbb{E}(\sup_{0\leq s\leq h} e^{-\alpha\xi_s})}.$$

Finally, Lemma B.1 leads to

$$\theta = h \frac{\mathbb{E}(\sup_{0\leq s\leq 1} e^{-\alpha\xi_s} - \sup_{s\geq 1} e^{-\alpha\xi_s})^+}{\mathbb{E}(\sup_{0\leq s\leq h} e^{-\alpha\xi_s})}, \tag{B.8}$$

which proves (a).

(b) We study the point process behavior of $(I_k)_{k\in\mathbb{N}}$ as in (a) by proving that the assumptions of Davis and Hsing (1995), Theorem 2.7, are satisfied. The finite dimensional distributions of $(I_k)_{k\in\mathbb{N}}$ are multivariate regularly varying by Theorem 3.3(c). At the end we will show that condition (B.3) for $(I_k)_{k\in\mathbb{N}}$ is satisfied. We define

$$G_k = \sqrt{V_{(k-1)h-}} e^{\xi_{(k-1)h-}/2} \int_{(k-1)h}^{kh} e^{-\xi_{t-}/2} \, dL_t,$$

$$A_k^* = e^{-(\xi_{(k-1)h-} - \xi_{h-})} \left( e^{\xi_{(k-1)h-}/2} \int_{(k-1)h}^{kh} e^{-\xi_{t-}/2} \, dL_t \right)^2,$$

$$B_k^* = e^{-\xi_{(k-1)h-}} \left( \int_h^{(k-1)h-} e^{\xi_{s-}} \, d\eta_s \right) \left( e^{\xi_{(k-1)h-}/2} \int_{(k-1)h}^{kh} e^{-\xi_{t-}/2} \, dL_t \right)^2, \qquad k\geq 2,$$

and $G_1^* = \max\{V_h, G_1^2, I_1^2\}$. Then $G_k^2 \leq A_k^* G_1^* + B_k^*$. Next, let $0 < \epsilon < x$, then

$$\mathbb{P}\left( \bigvee_{l\leq k\leq r_n} |I_k| > a_n^{1/2} x, |I_1| > a_n^{1/2} x \right)$$

$$\leq \sum_{l\leq k\leq r_n} \mathbb{P}(G_k^2 > a_n(x-\epsilon)^2, I_1^2 > a_n x^2) + \sum_{l\leq k\leq r_n} \mathbb{P}((I_k - G_k)^2 > a_n \epsilon^2, I_1^2 > a_n x^2)$$

$$=: I_{(n)} + II_{(n)}.$$

First, we investigate $I_{(n)}$. Note that $G_1^*$ is independent of $(A_k^*, B_k^*)$ for $k \geq 2$, and $G_1^*$ is regularly varying with index $\alpha$ and tail behavior

$$\mathbb{P}(G_1^* > x) \sim \mathbb{E}\left( \max\left\{ e^{-\alpha\xi_h}, \left| \int_0^h e^{-\xi_{t-}/2} \, dL_t \right|^{2\alpha} \right\} \right) \mathbb{P}(V_0 > x) \qquad \text{as } x \to \infty.$$

This is a conclusion of

$$G_1^* = \max\left\{ e^{-\xi_h}\left( V_0 + \int_0^h e^{\xi_{t-}} \, d\eta_t \right), \left( \int_0^h e^{-\xi_{t-}/2} \, dL_t \right)^2 V_0, \ldots, \right.$$

$$\left. \left( \int_0^h \sqrt{e^{-\xi_{t-}}\left( V_0 + \int_0^{t-} e^{\xi_{s-}} \, d\eta_s \right)} \, dL_t \right)^2 \right\}$$

and similar arguments as in Theorem 3.3(c). Furthermore,

$$\mathbb{E}(A_k^{*v}) = \mathbb{E}\left|\int_0^h e^{-\xi_{t-}/2}\,dL_t\right|^{2v} e^{\Psi_\xi(v)(k-2)h} \qquad \text{for } v \le d, k \ge 2.$$

Thus,

$$\mathbb{P}\left(\bigvee_{l \le k \le r_n} G_k^2 > a_n(x-\epsilon)^2, I_1^2 > a_n x^2\right)$$
$$\le \sum_{k=l}^{r_n} [\mathbb{P}(B_k^* > a_n(x-\epsilon)^2/2)\mathbb{P}(G_1^* > a_n x^2) + \mathbb{P}(A_k^* G_1^* > a_n(x-\epsilon)^2/2, G_1^* > a_n x^2)].$$

The remainder of the proof is as in (a) and we obtain

$$\lim_{l\to\infty}\lim_{n\to\infty}\frac{I_{(n)}}{\mathbb{P}(I_1^2 > a_n x^2)} = 0.$$

Next, we study $II_{(n)}$. By Markov's inequality we have for $\widetilde{d} = \max\{1,d\}$,

$$\mathbb{P}((I_k - G_k)^2 > a_n \epsilon^2, I_1^2 > a_n x^2) \le (a_n \epsilon^2)^{-\widetilde{d}}\mathbb{E}((I_k - G_k)^{2\widetilde{d}}\mathbf{1}_{\{I_1^2 > a_n x^2\}}).$$

Then similar computations as in (A.6) lead to

$$\mathbb{P}((I_k - G_k)^2 > a_n \epsilon^2, I_1^2 > a_n x^2)$$
$$\le K_1(a_n\epsilon^2)^{-\widetilde{d}}\mathbb{E}\left(\int_{(k-1)h}^{kh} e^{-\xi_{t-}}\left(\int_{(k-1)h}^{t-} e^{\xi_{s-}}\,d\eta_s\right)d[L,L]_t\right)^{\widetilde{d}}\mathbb{P}(I_1^2 > a_n x^2)$$
$$\le K_2(a_n\epsilon^2)^{-\widetilde{d}}\mathbb{E}\left|\int_{(k-1)h}^{kh} e^{-\xi_{t-}/2}\left(\int_{(k-1)h}^{t-} e^{\xi_{s-}}\,d\eta_s\right)^{1/2}dL_t\right|^{2\widetilde{d}}\mathbb{P}(I_1^2 > a_n x^2).$$

Thus, also

$$\lim_{l\to\infty}\lim_{n\to\infty}\frac{II_{(n)}}{\mathbb{P}(I_1^2 > a_n x^2)} = 0.$$

Hence, Theorem 2.7 in Davis and Hsing (1995) proves the statement. $\square$

**Proof of Lemma 4.5.** Let $(Z_k)_{k\in\mathbb{N}}$ be the jump sizes of $L$. We define $\widetilde{Z}_k = \log(1 + \lambda e^c Z_k^2)$ for $k \in \mathbb{N}$ and $\widetilde{Z}$ as a random variable with $\widetilde{Z} \stackrel{d}{=} \widetilde{Z}_1$. Then $\xi_t = ct - \sum_{k=1}^{N_t}\widetilde{Z}_k$, where $(N_t)$ is a Poisson process with jumps $(\Gamma_k)_{k\in\mathbb{N}}$, and by (4.1) in Klüppelberg *et al.* (2004),

$$\Psi_\xi(s) = -(\mu + sc) + \mu\mathbb{E}(e^{s\widetilde{Z}}) \qquad \text{for } s \le \alpha.$$

Thus, we obtain with $\Psi_\xi(\alpha) = 0$ that

$$\mathbb{E}(e^{\alpha \widetilde{Z}}) = \frac{\mu + \alpha c}{\mu} \quad \text{and} \quad \mathbb{E}(e^{-(\mu/c)\widetilde{Z}}) = \frac{1}{\mu}\Psi_\xi\left(-\frac{\mu}{c}\right). \tag{B.9}$$

Let $(\widetilde{\xi}_t)_{t\geq 0}$ be a Lévy process independent of $\xi$ and identically distributed as $\xi$. We write $\underline{\widetilde{\xi}}_t = \inf_{0\leq s\leq t} \widetilde{\xi}_s$, and $e_t$ for an exponentially distributed random variable with mean $1/t$ for $t > 0$, which is independent of $\xi$ and $\widetilde{\xi}$.

First, we investigate the right-hand side of (4.1). The following equality holds:

$$\mathbb{E}\left(1 - \sup_{s\geq \Gamma_1} e^{-\alpha \xi_s}\right)^+ = \mathbb{E}((1 - e^{-\alpha \xi_{\Gamma_1}} e^{-\alpha \underline{\widetilde{\xi}}_\infty})^+ \mathbf{1}_{\{\xi_{\Gamma_1}>0\}}).$$

Kyprianou (2006), Exercise 1.8(iii), says that

$$\mathbb{P}(\xi_{\Gamma_1} \mathbf{1}_{\{\xi_{\Gamma_1}>0\}} > x) = \mathbb{E}(e^{-(\mu/c)\widetilde{Z}})\mathbb{P}(e_{\mu/c} > x) \qquad \text{for } x > 0.$$

If we use (B.9), this leads to

$$\begin{aligned}\mathbb{E}\left(1 - \sup_{s\geq \Gamma_1} e^{-\alpha \xi_s}\right)^+ &= \mathbb{E}(e^{-(\mu/c)\widetilde{Z}})\mathbb{E}(1 - e^{-\alpha e_{\mu/c}} e^{-\alpha \underline{\widetilde{\xi}}_\infty})^+ \\ &= \frac{1}{\mu}\Psi_\xi\left(-\frac{\mu}{c}\right)\mathbb{E}(1 - e^{-\alpha e_{\mu/c}} e^{-\alpha \underline{\widetilde{\xi}}_\infty})^+.\end{aligned} \tag{B.10}$$

Next, we look at the left-hand side of (4.1). Let $q > 0$. By Lemma B.1 we have

$$\begin{aligned}\mathbb{E}\left(\sup_{s\leq 1} e^{-\alpha \xi_s} - \sup_{s\geq 1} e^{-\alpha \xi_s}\right)^+ &= q\mathbb{E}\left(\sup_{s\leq e_q} e^{-\alpha \xi_s} - \sup_{s\geq e_q} e^{-\alpha \xi_s}\right)^+ \\ &= q\mathbb{E}(e^{-\alpha \underline{\xi}_{e_q}}(1 - e^{-\alpha(\xi_{e_q}-\underline{\xi}_{e_q})} e^{-\alpha \underline{\widetilde{\xi}}_\infty}))^+.\end{aligned}$$

By the Wiener–Hopf decomposition (cf. Kyprianou (2006), Theorem 6.16), $\underline{\xi}_{e_q}$ and $\xi_{e_q} - \underline{\xi}_{e_q}$ are independent such that

$$\mathbb{E}\left(\sup_{s\leq 1} e^{-\alpha \xi_s} - \sup_{s\geq 1} e^{-\alpha \xi_s}\right)^+ = q\mathbb{E}(e^{-\alpha \underline{\xi}_{e_q}})\mathbb{E}(1 - e^{-\alpha(\xi_{e_q}-\underline{\xi}_{e_q})} e^{-\alpha \underline{\widetilde{\xi}}_\infty})^+.$$

Then (8.2) in Kyprianou (2006) results in

$$\begin{aligned}\mathbb{E}\left(\sup_{s\leq 1} e^{-\alpha \xi_s} - \sup_{s\geq 1} e^{-\alpha \xi_s}\right)^+ &= \Psi_\xi\left(-\frac{\mu}{c}\right)\frac{\mu + \alpha c}{\mu}\mathbb{E}(1 - e^{-\alpha e_{\mu/c}} e^{-\alpha \underline{\widetilde{\xi}}_\infty})^+ \\ &= \Psi_\xi\left(-\frac{\mu}{c}\right)(\mathbb{E}(e^{-\alpha c\Gamma_1}))^{-1}\mathbb{E}(1 - e^{-\alpha e_{\mu/c}} e^{-\alpha \underline{\widetilde{\xi}}_\infty})^+.\end{aligned} \tag{B.11}$$

The comparison of (B.10) and (B.11) gives the proof. □

We need the following lemma for the investigation of the convergence of the sample autocovariances.

**Lemma B.2.** *Let $(V_t)_{t\geq 0}$ be a genOU process satisfying (A) and (B), and define for $h > 0$, $\mathbf{V}_k = (V_{kh}, \ldots, V_{(k+m)h})$, $k \in \mathbb{N}_0$. Let $0 < a_n \uparrow \infty$ be a sequence of constants such that*

$$\lim_{n\to\infty} n\mathbb{P}(|\mathbf{V}_0|_\infty > a_n x) = x^{-\alpha} \qquad \text{for } x > 0.$$

*Suppose $(\mathbf{V}_k)_{k\in\mathbb{N}}$ satisfies $\mathcal{A}(a_n)$. Then*

$$\kappa_n := \sum_{k=1}^n \varepsilon_{a_n^{-1}\mathbf{V}_k} \stackrel{n\to\infty}{\Longrightarrow} \sum_{k=1}^\infty \sum_{j=0}^\infty \varepsilon_{\mathbf{Q}_{kj}P_k} =: \kappa, \tag{B.12}$$

*where $\sum_{k=1}^\infty \varepsilon_{P_k}$ is $\mathrm{PRM}(\vartheta)$ with*

$$\vartheta(\mathrm{d}x) = \alpha \frac{\mathbb{E}(\bigvee_{j=0}^m \mathrm{e}^{-\alpha\xi_{jh}} - \bigvee_{j=1}^\infty \mathrm{e}^{-\alpha\xi_{jh}})^+}{\mathbb{E}(\bigvee_{j=0}^m \mathrm{e}^{-\alpha\xi_{jh}})} x^{-\alpha-1} \mathbf{1}_{(0,\infty)}(x)\,\mathrm{d}x.$$

*Furthermore, $\sum_{j=0}^\infty \varepsilon_{\mathbf{Q}_{kj}}$ for $k \in \mathbb{N}$ are i.i.d. point processes independent of $\sum_{k=1}^\infty \varepsilon_{P_k}$ with $0 \leq |\mathbf{Q}_{kj}|_\infty \leq 1$, for each $k$ exactly one $|\mathbf{Q}_{kj}|_\infty$ is equal to 1, and $\mathbb{P}(|\mathbf{Q}_{kj}|_\infty = 0) < 1$ for $j \in \mathbb{N}$. The sequence $(\mathbf{Q}_{kj})_{j\in\mathbb{N}_0}$ is a.s. unique.*

One can either interpret $(\mathbf{V}_k)_{k\in\mathbb{N}}$ as a multivariate stochastic recurrence equation and apply Theorem 2.10 of Basrak *et al.* (2002) to obtain the proof of Lemma B.2, or one can proceed as in the proof of Theorem 4.1.

**Proof of Theorem 4.8.** Without loss of generality we can assume that $h = 1$.

(a, c) are conclusions of Lemma B.2, and Davis and Mikosch (1998), Theorem 3.5, and arguments presented on page 2069 there.

(b) The proof is similar to the proof in Mikosch and Stărică (2000), page 1440 ff., so that we present only a sketch of it. Let $\mathbf{x}_k = (x_k^{(0)}, \ldots, x_k^{(m)}) \in \overline{\mathbb{R}}^{m+1} \setminus \{\mathbf{0}\}$. We define the mappings $T_{j,\epsilon}: \mathcal{M} \to \overline{\mathbb{R}}$ by

$$T_{0,\epsilon}\left(\sum_{k=1}^\infty n_k \varepsilon_{\mathbf{x}_k}\right) = \sum_{k=1}^\infty n_k (x_k^{(0)})^2 \mathbf{1}_{\{|x_k^{(0)}|>\epsilon\}},$$

$$T_{1,\epsilon}\left(\sum_{k=1}^\infty n_k \varepsilon_{\mathbf{x}_k}\right) = \sum_{k=1}^\infty n_k (x_k^{(1)})^2 \mathbf{1}_{\{|x_k^{(0)}|>\epsilon\}},$$

$$T_{j,\epsilon}\left(\sum_{k=1}^\infty n_k \varepsilon_{\mathbf{x}_k}\right) = \sum_{k=1}^\infty n_k x_k^{(0)} x_k^{(j-1)} \mathbf{1}_{\{|x_k^{(0)}|>\epsilon\}}, \qquad j \geq 2.$$

Furthermore, we define $A_k = A_{k+1}^k$ and $B_k = B_{k+1}^k$ so that $V_{k+1} = A_k V_k + B_k$ for $k \in \mathbb{N}$, and $(A_k, B_k)_{k \in \mathbb{N}}$ is an i.i.d. sequence. First, we derive the asymptotic behavior of the sample variance

$$n^{1-2/\alpha}(\gamma_{n,V}(0) - \gamma_V(0))$$
$$= n^{-2/\alpha} \sum_{k=1}^n (V_{k+1}^2 - \mathbb{E}(V_k^2)) + o_p(1)$$
$$= n^{-2/\alpha} \sum_{k=1}^n V_k^2(A_k^2 - \mathbb{E}(A_k^2)) + \mathbb{E}(A_1^2) n^{-2/\alpha} \sum_{k=1}^n (V_k^2 - \mathbb{E}(V_k^2))$$
$$+ 2n^{-2/\alpha} \sum_{k=1}^n (A_k B_k V_k - \mathbb{E}(A_k B_k V_k)) + n^{-2/\alpha} \sum_{k=1}^n (B_k^2 - \mathbb{E}(B_k^2)) + o_p(1).$$

Using the central limit theorem (CLT) for exponentially $\alpha$-mixing sequences (cf. Ibragimov and Linnik (1971), Theorem 18.5.3) (where we require $d > 4$ such that by Hölder's inequality $\mathbb{E}(A_k^2 B_k^2) \leq (\mathbb{E}(A_k^4))^{1/2}(\mathbb{E}(B_k^4))^{1/2} < \infty$), we obtain

$$(1 - \mathbb{E}(A_1^2))n^{1-2/\alpha}(\gamma_{n,V}(0) - \gamma_V(0))$$
$$= n^{-2/\alpha} \sum_{k=1}^n (V_k^2(A_k^2 - \mathbb{E}(A_k^2)))\mathbf{1}_{\{V_k > n^{1/\alpha}\epsilon\}} \quad \text{(B.13)}$$
$$+ n^{-2/\alpha} \sum_{k=1}^n (V_k^2(A_k^2 - \mathbb{E}(A_k^2)))\mathbf{1}_{\{V_k \leq n^{1/\alpha}\epsilon\}} + o_p(1)$$
$$=: I_{\epsilon,n}^{(1)} + II_{\epsilon,n}^{(1)} + o_p(1).$$

We proceed with the investigation of the behavior of $II_{\epsilon,n}^{(1)}$, which is the sum of uncorrelated random variables. Hence, according to Karamata's theorem (see Feller (1971), VIII.9, Theorem 1) for $n \to \infty$ holds

$$\text{Var}(II_{\epsilon,n}^{(1)}) = n^{-4/\alpha} n \mathbb{E}((A_k^2 - \mathbb{E}(A_k^2))^2)\mathbb{E}(V_k^4 \mathbf{1}_{\{V_k \leq n^{1/\alpha}\epsilon\}}) \sim K\epsilon^{4-\alpha} \xrightarrow{\epsilon \downarrow 0} 0. \quad \text{(B.14)}$$

Let $\kappa$ and $\kappa_n$ be as in Lemma B.2. We denote by $(S_0^*, \ldots, S_m^*)$ the weak limit of

$$(T_{1,\epsilon}\kappa_n - \mathbb{E}(A_1^2)T_{0,\epsilon}\kappa_n, T_{2,\epsilon}\kappa_n - \mathbb{E}(A_1)T_{1,\epsilon}\kappa_n, \ldots, T_{m+1,\epsilon}\kappa_n - \mathbb{E}(A_1)T_{m,\epsilon}\kappa_n) =: T_\epsilon \kappa_n$$

as first $n \to \infty$ and then $\epsilon \downarrow 0$, which exists due to Lemma B.2, an extended version of (B.14), $\mathbb{E}(T_\epsilon \kappa) = 0$ and the arguments presented in Davis and Hsing (1995), proof of Theorem 3.1(ii), page 897 (cf. Mikosch and Stărică (2000), page 1441), that is,

$$T_\epsilon \kappa_n \xRightarrow{n \to \infty, \epsilon \downarrow 0} (S_0^*, \ldots, S_m^*). \quad \text{(B.15)}$$

For summand $I_{\epsilon,n}^{(1)}$ we obtain

$$\begin{aligned}
I_{\epsilon,n}^{(1)} =\ & n^{-2/\alpha}\sum_{k=1}^{n}((A_kV_k+B_k)^2 - \mathbb{E}(A_k^2)V_k^2)\mathbf{1}_{\{V_k>n^{1/\alpha}\epsilon\}} \\
& - n^{-2/\alpha}\sum_{k=1}^{n}(B_k^2\mathbf{1}_{\{V_k>n^{1/\alpha}\epsilon\}} - \mathbb{E}(B_k^2\mathbf{1}_{\{V_k>n^{1/\alpha}\epsilon\}})) \\
& - 2n^{-2/\alpha}\sum_{k=1}^{n}(A_kB_kV_k\mathbf{1}_{\{V_k>n^{1/\alpha}\epsilon\}} - \mathbb{E}(A_kB_kV_k\mathbf{1}_{\{V_k>n^{1/\alpha}\epsilon\}})) \\
& - n^{-2/\alpha}\mathbb{E}(B_1^2)n\mathbb{P}(V_0>n^{1/\alpha}\epsilon) - 2n^{-2/\alpha}\mathbb{E}(A_1B_1)n\mathbb{E}(V_1\mathbf{1}_{\{V_1>n^{1/\alpha}\epsilon\}}).
\end{aligned} \qquad \text{(B.16)}$$

A consequence of the CLT and the regular variation of $\mathbb{E}(V_1\mathbf{1}_{\{V_1>x\}})$ with index $\alpha-1$ as $x\to\infty$ is that

$$\begin{aligned}
I_{\epsilon,n}^{(1)} &= n^{-2/\alpha}\sum_{k=1}^{n}(V_{k+1}^2 - \mathbb{E}(A_k^2)V_k^2)\mathbf{1}_{\{V_k>n^{1/\alpha}\epsilon\}} + \mathrm{o}_p(1) \\
&= T_{1,\epsilon}\kappa_n - \mathbb{E}(A_1^2)T_{0,\epsilon}\kappa_n + \mathrm{o}_p(1) \\
&\stackrel{n\to\infty}{\Longrightarrow} T_{1,\epsilon}\kappa - \mathbb{E}(A_1^2)T_{0,\epsilon}\kappa \stackrel{\epsilon\downarrow 0}{\Longrightarrow} S_0^*.
\end{aligned} \qquad \text{(B.17)}$$

Equations (B.13), (B.14), (B.17) and Billingsley (1999), Theorem 3.1, lead to

$$n^{1-2/\alpha}(\gamma_{n,V}(0) - \gamma_V(0)) \stackrel{n\to\infty}{\Longrightarrow} (1 - \mathbb{E}(A_1^2))^{-1}S_0^* =: S_0. \qquad \text{(B.18)}$$

In the same way it is possible to extend the result to sample autocovariance functions of higher orders, where

$$n^{1-2/\alpha}(\gamma_{n,V}(l) - \gamma_V(l)) \stackrel{n\to\infty}{\Longrightarrow} S_l^* + \mathbb{E}(A_1)S_{l-1} =: S_l \qquad \text{for } l\geq 1. \qquad \text{(B.19)}$$

This results in (4.5).

We obtain the asymptotic behavior of the sample autocorrelation function from the behavior of the sample autocovariance function and the continuous mapping theorem as in Davis and Mikosch (1998), page 2061. $\square$

## Acknowledgements

Parts of the paper were written while the author was visiting the Department of Operations Research and Industrial Engineering at Cornell University and the School of Finance and Applied Statistics at the Australian National University. She takes pleasure in thanking them for their hospitality. Discussions with Alexander Lindner have been very fruitful and stimulating. Financial support from the Deutsche Forschungsgemeinschaft through a research grant is gratefully acknowledged.